\documentclass[a4paper,11pt,dvipdfmx]{article}

\makeatletter
\usepackage{amsmath,amsmath,amssymb,amsthm}

\usepackage{bm}
\usepackage{natbib}
\usepackage{titling}

\usepackage[dvipdfmx]{graphicx,color}
\usepackage{multirow}
\usepackage{xr}

\usepackage[plain,noend]{algorithm2e}
\newtheorem{theorem}{Theorem}
\newtheorem{lemma}{Lemma}
\newtheorem{definition}{Definition}
\newtheorem{corollary}{Corollary}
\newtheorem{algo}{Algorithm}
\newcommand{\sk}[1]{{\color{black}#1}}
\def\Real{\mathbb{R}}

\begin{document}
\title{\textbf{\Large M\"obius transformation and a Cauchy family on the sphere} \vspace{0.2cm} \\
}
\author{{\scshape Shogo Kato\thanks{\footnotesize {\it Address for correspondence}: Shogo Kato, Institute of Statistical Mathematics, 10-3 Midori-cho, Tachikawa, Tokyo 190-8562, Japan. {\tt E-mail:\ skato@ism.ac.jp}} \(^{,a}\) \ and Peter\ McCullagh\,\(^{b}\)} \vspace{0.5cm}\\
\textit{\(^{a}\) Institute of Statistical Mathematics, Japan} \\  
\textit{\(^{b}\) University of Chicago, USA}
}
\date{}
\maketitle

\begin{abstract}
We present some properties of a Cauchy family of distributions on the sphere, which
is a spherical extension of the wrapped Cauchy family on the circle.
The spherical Cauchy family is closed under the M\"obius transformation on the sphere and there is a similar induced transformation on the parameter space.
Stereographic projection transforms the the spherical \sk{Cauchy} family into a multivariate $t$-family \sk{with a certain degree of freedom} on Euclidean space.
Many tractable properties of the spherical Cauchy are derived using the M\"obius transformation and stereographic projection.
A method of moments estimator and an asymptotically efficient estimator are expressed in closed form.
The maximum likelihood estimation is also straightforward.\vspace{0.4cm}

\end{abstract}

\section{Introduction}

This paper discusses a family of distributions on the sphere $S^d \subset \Real^{d+1}$ with probability density function
\begin{equation}
f(y;\mu,\rho) = \frac{\Gamma \{(d+1)/2\}}{2 \pi^{(d+1)/2}} \left( \frac{ 1- \rho^2 }{1+\rho^2-2 \rho \mu^T x} \right)^d, \quad y \in S^d, \label{eq:sphere_c}
\end{equation}
with respect to surface area, where $\mu \in S^d$ is the location parameter, $\rho \in [0,1)$ is the concentration parameter, 
and $S^d = \{ x \in \mathbb{R}^{d+1} \, ; \, \|x\|=1 \}$ denotes the unit sphere in $\mathbb{R}^{d+1}$.
The circular case ($d=1$) is well-known as the wrapped Cauchy or circular Cauchy family; 
see, e.g., \citet{ken88} and \citet{mcc96}.
In this paper, the distribution (\ref{eq:sphere_c}) is called the Cauchy distribution on the sphere or the spherical Cauchy distribution.

\cite{mcc96} showed that the wrapped Cauchy family is closed under conformal maps preserving the unit circle \sk{which are called the M\"obius transformations on the unit circle}, and that there is a similar induced transformation on the parameter space.
Related results about the Cauchy family on the real line and on the Euclidean space have been given by \citet{mcc92} and \citet{let}, respectively.
To our knowledge, however, there has been no literature about the association between the M\"obius transformation and the spherical Cauchy family (\ref{eq:sphere_c}).
Since there have been various statistical applications of the wrapped Cauchy family and/or the M\"obius transformation in directional statistics \citep{mcc96,dow02,dow03,jon,kat10a,kat10b,kat15,ues}, it is potentially useful to consider the Cauchy family on the sphere and its relationship with the M\"obius transformation.


This paper presents some properties of the Cauchy family on the sphere, especially, those related to M\"obius transformation.
The spherical Cauchy family is closed under the M\"obius transformation on the sphere,
and the transformed parameter is given by the extended M\"obius transformation on $\Real^{d+1}$.
The statistical benefits of this property include: 
(i)~an efficient algorithm for random variate generation;
(ii)~a simple pivotal statistic for parametric inference;
(iii)~straightforward calculation of probabilities of a surface region;
(iv)~closed form expression for maximum likelihood estimator for $n \leq 3$; and 
(v)~straightforward calculation of the Fisher information matrix.
A method of moments estimator can be expressed in simple form.
A simple algorithm for maximum likelihood estimation is available.
The likelihood for the spherical Cauchy is equivalent to that for the $t$-family with a certain degree of freedom which is related to the spherical Cauchy via stereographic projection.
An asymptotically efficient estimator is presented which our simulation study suggests outperforms the method of moments estimator and the maximum likelihood estimator in certain settings.
Comparing the densities of the spherical Cauchy and von Mises--Fisher, the spherical Cauchy density takes greater values around the mode and antimode and smaller values in the other area of the sphere.
The advantages of the spherical Cauchy over the von Mises--Fisher in terms of properties include the closure under the M\"obius transformation and the related properties, while the von Mises--Fisher \sk{compares} favourably with the spherical Cauchy in terms of its membership in the exponential family, straightforward maximum likelihood estimation and well-developed theory of hypothesis testing.

Throughout this paper, $d$ is a positive integer.
We let $\mathbb{R}^d$ and $\overline{\mathbb{R}}^d$ denote the $d$-dimensional Euclidean space and extended Euclidean space $\mathbb{R}^d \cup \{\infty\}$, respectively.
Suppose that $\| \cdot \|$ is the Euclidean norm and that $S^d$ is the $d$-dimensional unit sphere in \sk{$\mathbb{R}^{d+1}$}, namely, $S^d = \{ x \in \mathbb{R}^{d+1} \, ; \, \|x\|=1 \}$.
Let $D^d$ and $\overline D^d$ denote the open and closed unit balls in $\Real^{d+1}$, so that
$S^d = \overline D^d \setminus D^d$.
The set of all $(d+1) \times (d+1)$ rotation matrices is denoted by $SO(d+1)$.
The $(d+1)$-dimensional unit vector whose $j$th element equals one is $e_j$.
The $(d+1) \times (d+1)$ identity matrix is denoted by $I$.

Proofs and further details can be found in the Supplementary Material.

\section{M\"obius transformation and a Cauchy family on $S^d$} \label{sec:cauchy_sphere}

\subsection{M\"obius transformation on $S^d$} \label{sec:mobius_sphere}

\sk{The goal of this section is to discuss the M\"obius transformation $S^d \to S^d$,
and to investigate its association with the Cauchy family~(\ref{eq:sphere_c}).
The first step to achieve this is to consider the following function
\begin{equation}
\tilde{\cal M}_{R,\sk{\psi}}(y) =  R \left\{ \frac{1- \| \sk{\psi} \|^2}{\|y+\sk{\psi} \|^2} (y+\sk{\psi}) + \sk{\psi} \right\}, \quad y \in S^d, \label{eq:mobius_s}
\end{equation}
where $\sk{\psi} \in \mathbb{R}^{d+1} \setminus S^d$ and $R \in SO(d+1)$.
The transformation (\ref{eq:mobius_s}) maps the unit sphere onto itself:
it is called the M\"obius transformation on the sphere.

The transformation (\ref{eq:mobius_s}) with $R=I$ can be interpreted as follows.
To any interior point $\|\sk{\psi}\| < 1$ there corresponds an inversion $S^d \to S^d$ that sends each point
$y\in S^d$ to an antipodal point $\tilde y$ by projection through~$\sk{\psi}$:
\begin{equation}\label{eq:inversion}
\tilde y - \sk{\psi} = - \frac{1 - \|\sk{\psi}\|^2}{\|y - \sk{\psi}\|^2} ( y - \sk{\psi}).
\end{equation}
The vectors $y-\sk{\psi}$ and $\tilde y-\sk{\psi}$ are co-linear, but opposite in direction, so $\sk{\psi}$ lies on the line segment $(y, \tilde y)$.
The product of the lengths is constant  $\|y-\sk{\psi}\| \times \|\tilde y - \sk{\psi}\| = 1 - \|\sk{\psi}\|^2$.
It follows from the intersecting chords theorem that $\tilde y \in S^d$.
As is well known, the intersecting chords theorem applies also to chords that intersect outside the circle,
so the transformation extends to~$\|\sk{\psi}\| > 1$.
In either case, the transformation is an inversion because $(y, \tilde y)$ is a $\sk{\psi}$-antipodal pair,
and a second application of (\ref{eq:inversion}) returns the original point
i.e.,~$\tilde{\tilde y} = y$ for every $y\in S^d$.
The inversion (\ref{eq:inversion}) is said to be conformal because there is no local distortion of angles:
it is a property of conic sections that the image of a circular cap is a circular cap.}

A function equivalent to (\ref{eq:mobius_s}) with the restriction $R=I$, \sk{$\sk{\psi}=\sk{\psi}_1 e_1$} and $-1<\sk{\psi}_1<1$
occurs in \S 10 of \cite{mcc89}.

The two parameters $R$ and $\sk{\psi}$ have a clear interpretation.
The matrix $R$ works as a rotation parameter.
In order to discuss the interpretation of $\sk{\psi}$, assume, without loss of generality, that $R=I$.
If $\| \sk{\psi} \|<1$, $\sk{\psi}$ can be interpreted as a parameter vector that attracts the points on the sphere towards $\sk{\psi}/\| \sk{\psi} \|$, with the concentration of the points around $\sk{\psi}/\| \sk{\psi} \|$ increasing as $\|\sk{\psi} \|$ increases.
In particular, if $\sk{\psi}=0$, then \sk{$\tilde{\cal M}_{I,\sk{\psi}}$} reduces to the identify mapping.
As $\| \sk{\psi} \| \rightarrow 1$, \sk{$ \tilde{\cal M}_{I,\sk{\psi}}(y)$} $\rightarrow$ $\sk{\psi}/\| \sk{\psi} \|$ for any $y \neq -\sk{\psi}/ \| \sk{\psi} \|$.
\sk{The points $y=\sk{\psi}/\| \sk{\psi} \|$ and $y=-\sk{\psi}/\| \sk{\psi} \|$ are invariant under $\tilde{\cal M}_{I,\sk{\psi}}$, i.e., $\tilde{\cal M}_{I,\sk{\psi}}(\sk{\psi}/\| \sk{\psi} \|)=\sk{\psi}/ \| \sk{\psi} \|$ and $\tilde{\cal M}_{I,\sk{\psi}} (-\sk{\psi}/ \| \sk{\psi} \|)=-\sk{\psi}/ \| \sk{\psi} \|$ for any $\sk{\psi} \neq 0$.}
For the case of $\| \sk{\psi}\|>1$, the transformation \sk{$\tilde{\cal M}_{I,\sk{\psi}}$} consists of the two types of transformations, namely, the reflection in $y=c \sk{\psi}/\|\sk{\psi}\| \ (c \in \mathbb{R})$ and the transformation \sk{$\tilde{\cal M}_{I,\sk{\psi}/\|\sk{\psi} \|^2}$}.

\subsection{An extension of the M\"obius transformation on $S^d$} \label{sec:extended_mobius}


\begin{definition}
We define a function by
\begin{equation}
{\cal M}_{R,\sk{\psi}}(x) =  R \left\{ \frac{1- \| \sk{\psi} \|^2}{\|\tilde{x}+\sk{\psi} \|^2} (\tilde{x}+\sk{\psi}) + \sk{\psi} \right\}, \quad x \in \mathbb{R}^{d+1} \setminus \{ 0,- \sk{\psi} /\| \sk{\psi}\|^2 \}. \label{eq:mobius_s2}
\end{equation}
where $\tilde{x} = x/\|x\|^2$, $\sk{\psi} \in \mathbb{R}^{d+1} \sk{\setminus S^d}$, and $R \in SO(d+1)$.
Also, we define ${\cal M}_{R,\sk{\psi}}(0)=R \sk{\psi}$, ${\cal M}_{R,\sk{\psi}}(-\sk{\psi}/\|\sk{\psi}\|^2)=\infty$ and ${\cal M}_{R,\sk{\psi}}(\infty) = R \sk{\psi}/\|\sk{\psi}\|^2$.
\end{definition}

If we restrict the domains of $x$ \sk{to be $S^d$}, then ${\cal M}_{R,\sk{\psi}}$ reduces to \sk{the M\"obius transformation on the sphere ${\cal M}_{R,\sk{\psi}}$}.
The transformation \sk{${\cal M}_{R,\sk{\psi}}$} can also be expressed as
\begin{equation}
{\cal M}_{R,\sk{\psi}}(x) = R \, T_{\tilde{\sk{\psi}}} \left\{ \frac{1- \| \tilde{\sk{\psi}} \|^2}{\| x + \tilde{\sk{\psi}} \|^2} ( x+ \tilde{\sk{\psi}}) + \tilde{\sk{\psi}} \right\}, \quad x \in \mathbb{R}^{d+1} \setminus \{ - \tilde{\sk{\psi}} \}, \label{eq:mobius_s3}
\end{equation}
where $\tilde{\sk{\psi}}=\sk{\psi}/\|\sk{\psi}\|^2$ and \sk{$T_{\tilde{\sk{\psi}}}=2 \tilde{\sk{\psi}} \tilde{\sk{\psi}}^T /\| \tilde{\sk{\psi}} \|^2-I$.}
Throughout the paper the transformation (\ref{eq:mobius_s2}) is denoted by ${\cal M}_{R,\sk{\psi}}$. 

\begin{theorem} \label{thm:mobius_conformal}
The following hold for the transformation ${\cal M}_{R,\sk{\psi}}$:
\begin{enumerate}
\item[(i)] The transformation ${\cal M}_{R,\sk{\psi}}$ is a bijective conformal map which maps $\overline{\mathbb{R}}^{d+1}$ onto itself.
\item[(ii)] For any $\sk{\psi} \in \overline{\mathbb{R}}^{d+1} \setminus S^d$, the transformation ${\cal M}_{R,\sk{\psi}}$ maps the unit sphere $S^d$ onto itself.
\item[(iii)] If $\|\sk{\psi}\|<1$, then ${\cal M}_{R,\sk{\psi}}(D^{d+1}) = D^{d+1}$ and ${\cal M}_{R,\sk{\psi}}(\overline{\mathbb{R}}^{d+1} \setminus \overline{D}^{d+1} ) = \overline{\mathbb{R}}^{d+1} \setminus \overline{D}^{d+1} $.
\item[(iv)] If $\|\sk{\psi}\|>1$, then ${\cal M}_{R,\sk{\psi}}(D^{d+1}) = \overline{\mathbb{R}}^{d+1} \setminus \overline{D}^{d+1}$ and ${\cal M}_{R,\sk{\psi}}(\overline{\mathbb{R}}^{d+1} \setminus \overline{D}^{d+1}) =  D^{d+1} $.
\end{enumerate}
\end{theorem}

If $d=1$, the transformation ${\cal M}_{R,\sk{\psi}}$ is related to the M\"obius transformation on the complex plane which is of the form
\begin{equation}
{\cal M}_c(z) = \alpha_0 \frac{z+\alpha_1}{\overline{\alpha}_1 z +1}, \quad z \in \overline{\mathbb{C}},  \label{eq:mobius_plane}
\end{equation}
where $\alpha_0$ and $\alpha_1$ are complex numbers such that $|\alpha_0|=1$ and $|\alpha_1| \neq 1$.
The transformation \sk{${\cal M}_c$} is essentially the same as \sk{${\cal M}_{R,\sk{\psi}}$} with $d=1$ if the real and imaginary parts of \sk{${\cal M}_c(z)$} are identified as the first and second components of ${\cal M}_{R,\sk{\psi}}(x)$, respectively.
\sk{This fact can be easily confirmed by expressing (\ref{eq:mobius_plane}) as
$$
{\cal M}_c (z) = \alpha_0 \frac{\alpha_1^2}{|\alpha_1|^2} \overline{ \left\{ \frac{1-|\alpha_1|^{-2}}{|z+\alpha_1/|\alpha_1|^2|^2} \left( z + \frac{\alpha_1}{|\alpha_1|^2} \right) + \frac{\alpha_1}{|\alpha_1|^2} \right\} }.
$$}

After extremely tedious but straightforward calculation, it follows that the set of transformations (\ref{eq:mobius_s3}) has the following closure property.
\begin{lemma} \label{lem:closure_mobius_s}
For $\sk{\psi}_2 \neq -R_1 \sk{\psi}_1 $,
$$
{\cal M}_{R_2,\sk{\psi}_2} \{ {\cal M}_{R_1, \sk{\psi}_1} (x) \} = {\cal M}_{\check{R},\check{\sk{\psi}}} (x),  \quad x \in \overline{\mathbb{R}}^{d+1}, \label{eq:closure_mobius_s}
$$
where
$\check{R} = R_2 T_{\sk{\psi}_2} T_{\beta} R_1 T_{\sk{\psi}_1} T_{\check{\sk{\psi}}}$, $\check{\sk{\psi}} = {\cal M}_{I,\sk{\psi}_1}(\sk{R_1^T} \tilde{\sk{\psi}}_2) / \| {\cal M}_{I,\sk{\psi}_1}(\sk{R_1^T} \tilde{\sk{\psi}}_2) \|^2 = T_{\sk{\psi}_1} \sk{R_1^T} T_{\beta} T_{\sk{\psi}_2} {\cal M}_{I,\sk{\psi}_2} (R_1 \sk{\psi}_1) $, $\beta= R_1 \tilde{\sk{\psi}}_1+ \tilde{\sk{\psi}}_2$, and \sk{$\tilde{\sk{\psi}}_j = \sk{\psi}_j/\|\sk{\psi}_j\|^2 \ (j=1,2)$}. 
If $\sk{\psi}_2 = -R_1 \sk{\psi}_1$, then ${\cal M}_{R_2,\sk{\psi}_2} \{ {\cal M}_{R_1,\sk{\psi}_1}(x) \} = {\cal M}_{R_2 R_1,0}(x)$ for any $x \in \overline{\mathbb{R}}^{d+1} $.
\end{lemma}
Using this lemma, the following result can be immediately obtained.
\begin{theorem} \label{thm:group}
Let ${\cal G}$ be a set of the transformations $\{{\cal M}_{R,\sk{\psi}}\}$ with all possible combinations of $R \in SO(d+1)$ and $\sk{\psi} \in \mathbb{R}^{d+1} \sk{\setminus S^d}$, namely, ${\cal G} = \{ {\cal M}_{R,\sk{\psi}} \, ; \, R \in SO(d+1),\, \sk{\psi} \in \mathbb{R}^{d+1} \sk{\setminus S^d} \}$.
Then ${\cal G}$ forms a group under composition.
\end{theorem}

Therefore $\{{\cal M}_{R,\sk{\psi}}\}$ is a subgroup of the M\"obius group on $\overline{\mathbb{R}}^{d+1}$; see, e.g., \citet[\S 2]{iwa}.
\sk{It is clear from Theorem \ref{thm:group} that the set of M\"obius transformations on the sphere $\{\tilde{\cal M}_{R,\sk{\psi}} \}$ also forms a group under composition.}
The set of transformations $\{{\cal M}_{R,\sk{\psi}}\}$ is not an abelian group, implying that that ${\cal M}_{R_1,\sk{\psi}_1} \{ {\cal M}_{R_2,\sk{\psi}_2} (x) \} = {\cal M}_{R_2,\sk{\psi}_2} \{ {\cal M}_{R_1,\sk{\psi}_1} (x) \} $ does not hold in general.
However, for fixed $\mu \in S^d$, the subset of transformations \sk{$\{ {\cal M}_{I,\rho \mu} \, ; \, |\rho| \neq 1 \}$} forms an abelian group.
\sk{Similarly, an abelian group can be established for the set of the M\"obius transformations on the sphere $\{ \tilde{\cal M}_{I,\rho \mu} \, ; \, |\rho| \neq 1 \}$.}

\subsection{A Cauchy family on $S^d$} \label{subsec:cauchy_sphere}


The parameters $\mu$ and $\rho$ of the spherical Cauchy family (\ref{eq:sphere_c})  can be clearly interpreted.
The parameter $\mu$ controls the mode of the density.
The concentration of the distribution is regulated by $\rho$.
The greater the value of $\rho$, the greater the concentration of the density (\ref{eq:sphere_c}) around the mode.
In particular, when \sk{$\rho = 0$, the distribution (\ref{eq:sphere_c}) reduces to the uniform distribution on $S^d$.}
On the other hand, as $\rho$ tends to 1, the distribution converges to a point distribution with singularity at $y=\mu$.
\sk{There is a similar interpretation for the case $\|\phi\|>1$ because $f(y;\phi)=f(y;\phi/\|\phi \|^2)$.
See Fig.\ \ref{fig:comparison} given in \S \ref{sec:comparison} for some plots of the densities of the spherical Cauchy (\ref{eq:sphere_c}).}

In order to investigate the relationship between the spherical Cauchy family (\ref{eq:sphere_c}) and the set of transformations (\ref{eq:mobius_s2}), it is advantageous to write the parameters of the spherical Cauchy (\ref{eq:sphere_c}) as $\phi=\rho \mu$ and extend the parameter space to be $\overline{\mathbb{R}}^{d+1}$.
Specifically, we write the density of the spherical Cauchy as
\begin{equation}
f(y;\phi) = \frac{\Gamma \{(d+1)/2\}}{2 \pi^{(d+1)/2}} \left( \frac{ \left| 1- \|\phi \|^2 \right| }{\| y - \phi \|^2} \right)^d, \quad y \in S^d, \label{eq:sphere_c2}
\end{equation}
where $\phi \in \mathbb{R}^{d+1} \setminus S^d$.
For $\phi \in S^d$ we assume that the distribution is a point mass at $y=\phi$.
Also define that the density is uniform if $\phi=\infty$.
It can be seen that $f(y;\phi) = f(y;\phi/\|\phi\|^2)$ for any $\phi$.
Write $Y \sim C^*_{d}(\phi)$ if an $S^d$-valued random vector $Y$ has density (\ref{eq:sphere_c2}).

The following result can be readily established from Lemma \ref{lem:closure_mobius_s}.
\begin{theorem} \label{thm:closure_cauchy_s}
The following hold for the spherical Cauchy family (\ref{eq:sphere_c2}) and the transformation \sk{${\cal M}_{R,\psi}$}:
$$
Y \sim C_d^*(\phi) \quad \Longrightarrow \quad {\cal M}_{R,\psi} (Y) \sim C_d^* \left\{ {\cal M}_{R,\psi}(\phi) \right\}.
$$
\end{theorem}

If $d=1$, Theorem \ref{thm:closure_cauchy_s} is essentially the same as the result for the circular Cauchy or wrapped Cauchy family given in \cite{mcc96}.

There are some statistical applications of Theorem \ref{thm:closure_cauchy_s}.
For example, this theorem can be applied to propose an efficient algorithm to generate a random variate following the Cauchy family on $S^d$. 
\begin{corollary} \label{cor:random_variate}
If a random vector $U$ follows the uniform distribution on $S^d$, then 
 ${\cal M}_{I,\phi}(U)$ has the Cauchy distribution on the sphere $C^*_d(\phi)$.
\end{corollary}


In addition a pivotal statistics for $\phi$ and probabilities of a surface area under the density (\ref{eq:sphere_c2}) can be obtained as follows.
\begin{corollary}
Suppose $Y  \sim C^* (\phi)$.
Then ${\cal M}_{R,-\phi} (Y)$ is a pivotal statistics for $\phi$, where $R$ is any $(d+1)\times (d+1)$ matrix.
\end{corollary}

\begin{corollary} \label{cor:probabilities}
Let $f(y;\phi)$ denote the density (\ref{eq:sphere_c2}).
Assume $A \subset S^d$.
Then
$$
\int_A f(y;\phi) dy = \frac{Area \{{\cal M}_{I,-\phi} (A) \}}{Area(S^d)},
$$
where $Area(C)$ denotes the area of $C$ with respect to the surface measure.
\end{corollary}
The proofs of Corollaries \ref{cor:random_variate}--\ref{cor:probabilities} are straightforward from Theorem \ref{thm:closure_cauchy_s} and omitted.

\section{\sk{Extended stereographic projection}} \label{sec:inverse}

In this section we consider a transformation of the Cauchy family on the sphere (\ref{eq:sphere_c2}) via the stereographic projection.
The stereographic projection $S^d \rightarrow \overline{\mathbb{R}}^d$ is known to be
\begin{equation}
\tilde{\cal P}(y) = \frac{1}{1-y_{d+1}} (y_1,\ldots,y_d)^T, \quad y \in S^d \setminus \{e_{d+1}\}.  \label{eq:gis}
\end{equation}
Also assume $\tilde{\cal P}(e_{d+1}) = \infty$.
It is known that the stereographic projection (\ref{eq:gis}) maps the unit sphere $S^d$ onto $\overline{\mathbb{R}}^d$.
A geometrical interpretation of (\ref{eq:gis}) is that $\tilde{\cal P}(y)$ corresponds to the point at the intersection of the embedded Euclidean space $\overline{\mathbb{R}}^{d} \times \{0\}$ and the line connecting $y$ and the north pole $e_{d+1}$.

In order to discuss the transformation of the spherical Cauchy family (\ref{eq:sphere_c2}) via the stereographic projection (\ref{eq:gis}), we propose an extension of the complex number and an extended stereographic projection.
\begin{definition}
An extension of the complex number is defined by
$$
\theta = \mu + i \sigma,
$$ 
where $\mu \in \mathbb{R}^d$, $\sigma \in \mathbb{R}$ and $i$ is the imaginary number.
We write $\mu + i \sigma = \mu$ if $\sigma=0$.
\end{definition}
\begin{definition}
We define a new function on $\overline{\mathbb{R}}^{d+1}$ by
\begin{equation}
{\cal P}(x) = 2 \frac{(x_1,\ldots,x_d)^T}{\|x-e_{d+1}\|^2} + i \, \frac{1-\|x\|^2}{ \|x-e_{d+1}\|^2 }, \quad x \in \mathbb{R}^{d+1} \setminus \{e_{d+1}\}  \label{eq:ggs}
\end{equation}
Also, ${\cal P}(\infty) = -i $ and ${\cal P}(e_{d+1}) = \infty$.
\end{definition}

\begin{theorem} \label{thm:ggs}
The following hold for the function \sk{${\cal P}$}.
\begin{enumerate}
\item[(i)] The function \sk{${\cal P}$} is a bijective function which maps $\overline{\mathbb{R}}^{d+1}$ onto  $(\mathbb{R}^d + i \mathbb{R}) \cup \{\infty\}$.
\item[(ii)] The function \sk{${\cal P}$} reduces to \sk{$\tilde{\cal P}$} if $x \in S^d$.
\item[(iii)] If $\|x\|<1 \ (\|x\|>1)$, then the imaginary part of ${\cal P}(x)$ is positive (negative). 
\end{enumerate}
\end{theorem}


Theorem \ref{thm:ggs} implies that there exists the inverse function of (\ref{eq:ggs}) which is
$$
{\cal P}^{-1} (\theta) = \frac{2}{\| \mu \|^2 + (1+\sigma)^2} \left( \mu^T,  \frac{\| \mu\|^2+\sigma^2-1}{2} \right)^T,
$$
where $\theta = \mu+i \sigma \in (\mathbb{R}^d + i \mathbb{R}) \setminus \{-i\}$.
Also, suppose that ${\cal P}^{-1}(-i)=\infty$ and ${\cal P}^{-1}(\infty)=e_{d+1}$.
Then the following result is established.
\begin{theorem} \label{thm:closure_two_cauchy}
\sk{The following hold for the spherical Cauchy family (\ref{eq:sphere_c2}) and the extended stereographic projection ${\cal P}$:}
$$
Y \sim C^*_d (\phi) \quad \Longrightarrow \quad {\cal P}(Y) \sim C_d \left\{ {\cal P}(\phi) \right\}.
$$
Equivalently,
$$
X \sim C_d (\theta) \quad \Longrightarrow \quad {\cal P}^{-1}(X) \sim C^*_d \left\{ {\cal P}^{-1}(\theta) \right\}. \label{eq:two_c}
$$
Here $C_d(\theta)$ denotes a \sk{$d$-variate} $t$-distribution with $d$ degrees of freedom with density
\begin{equation}
f(x;\theta) = \frac{2^{d-1} \Gamma \{(d+1)/2\}}{\pi^{(d+1)/2}} \left( \frac{ |\sigma| }{ \sigma^2 + \| x - \mu \|^2} \right)^d, \quad x \in \overline{\mathbb{R}}^d,  \label{eq:euclid_c2}
\end{equation}
where $\theta=\mu + i \sigma$, $\mu \in \mathbb{R}^d$ and $\sigma \neq 0$.
For $\sigma=0$, we assume that the distribution (\ref{eq:euclid_c2}) is a point mass at $x=\mu$.
If $\theta=\infty$, then the distribution is assumed to be a point distribution with singularity at $x=\infty$.
\end{theorem}

This theorem and Theorem \ref{thm:closure_cauchy_s} imply that a random variate following the $t$-distribution with $d$ degrees of freedom $C_d(\theta)$ can be generated from the uniform distribution on $S_d$.


\section{Statistical inference}

\subsection{Method of moments estimation} \label{sec:mme}

Throughout this section we assume that $Y_1, \ldots, Y_n$ is a random sample from the multivariate Cauchy distribution on the sphere $C_d^*(\phi)$ with $\| \phi \| <1$.

\begin{theorem} \label{thm:moments_sc}
Let $Y$ have the spherical Cauchy $C_d^*(\phi)$ with $ \|\phi \| < 1$.
Then, for $\phi \neq 0$,
$$
E(Y) = \eta_{1,d} (\|\phi\|) \frac{\phi}{\| \phi \|}, \quad E(YY^T) = \frac{1}{d} \left[ \{ 1-\eta_{2,d} (\|\phi\|) \} I + \{ (d+1) \eta_{2,d} (\|\phi\|) -1 \} \frac{\phi \phi^T}{\|\phi\|^2}  \right],
$$
where
$$
\sk{\eta_{1,d} (\rho) } = \frac{1+\rho^2}{2 \rho} \left[ 1 - \frac{(1+\rho)^2}{1+\rho^2} \, F \left\{ 1 , \frac{d}{2}; d;  \frac{-4 \rho}{(1-\rho)^2} \right\} \right],
$$
$$
\sk{\eta_{2,d} (\rho) } = \frac{(1+\rho^2)^2}{4 \rho^2} \Biggl[ 1 - 2 \frac{(1+\rho)^2}{1+\rho^2} F \left\{ 1 ,\frac{d}{2}; d; \frac{-4 \rho}{(1-\rho)^2} \right\} + \frac{(1+\rho)^4}{(1+\rho^2)^2} F \left\{ 2 ,\frac{d}{2} ; d ; \frac{-4 \rho}{(1-\rho)^2} \right\} \Biggr],
$$
and $F$ denotes the hypergeometric series \citep[equation 9.111]{gra}.
If $\phi=0$, $E(Y)=0$ and $E(YY^T) = (d+1)^{-1} I$.
\end{theorem}
This theorem and Theorems \ref{thm:1st_moment} and \ref{thm:2nd_moment} of the Supplementary Material imply that $E(Y)$ and $E(YY^T)$ can be expressed in closed form without hypergeometric functions for any $d$.

A method of moments estimator of $\phi$ is obtained by equating the expectation of $Y$ and its sample analogue.
In other words the method of moments estimator is the solution of the equation
\begin{equation}
\eta_{1,d}(\| \phi \|) \frac{\phi}{\| \phi \|} = \overline{Y}, \label{eq:mom}
\end{equation}
where $\eta_{1,d} (\|\phi\|)$ is defined as in Theorem \ref{thm:moments_sc} and $ \overline{Y} = n^{-1} \sum_{j=1}^n Y_j$.
\sk{As is clear from Lemma \ref{lem:marginal} of Supplementary Material, it holds that $\eta_{1,d}(0) = 0$,  $\lim_{\| \phi \| \rightarrow 1} \eta_{1,d}(\| \phi \|) = 1$, and $\eta_{1,d}(\| \phi \|)$ is monotonically increasing with respect to $\|  \phi \|$.
This immediately leads to the following theorem.}
\begin{theorem}
The equation (\ref{eq:mom}) has the unique solution $\hat{\phi}_{MM}$ on the $(d+1)$-dimensional unit disc
\begin{equation}
\hat{\phi}_{MM} = \eta_{1,d}^{-1} \left( \| \overline{Y} \| \right) \frac{\overline{Y}}{\|\overline{Y}\|}, \label{eq:mme}
\end{equation}
where \sk{$\eta_{1,d}^{-1}(\rho)$ is the inverse of $\eta_{1,d} (\rho)$ for $0 \leq \rho <1$}.
\end{theorem}
Since $\eta_{1,d}$ is monotonically increasing, the method of moments estimate $\hat{\phi}_{MM}$ can be estimated numerically via usual optimization algorithms.

\begin{theorem} \label{thm:mme}
Let $\hat{\phi}_{MM}$ be the method of moments estimator (\ref{eq:mme}).
Then $ \sqrt{n} (\hat{\phi}_{MM} - \phi) $ tends in distribution to $N(0, \Lambda \Sigma \Lambda )$ as $n \rightarrow \infty$, where 
$$
\Sigma = d^{-1} \left[ \{ 1 - \eta_{2,d} (\|\phi\|) \} I + \{ (d+1) \eta_{2,d} (\| \phi \|) -1 -d \eta_{1,d}^2(\| \phi \|) \} \frac{\phi \phi^T}{\|\phi\|^2} \right],
$$
$$
\Lambda = {\eta_{1,d}^{-1}}' \{ \eta_{1,d} (\|\phi\|) \} \frac{\phi^T \phi}{\|\phi\|^2} + \frac{\|\phi\|}{|\eta_{1,d} (\|\phi\|)|} \left( I - \frac{\phi^T \phi}{ \|\phi\|^2} \right),
$$
$$
{\eta_{1,d}^{-1}}' \{ \eta_{1,d} (\|\phi\|) \} = \frac{d+1}{2d} \frac{(1-\|\phi\|)^3}{1+\|\phi\|} \left\{ F \left( 2 , \frac{d}{2}+1 ; d+2 ; - \frac{4\|\phi\|}{(1-\|\phi\|)^2} \right) \right\}^{-1}. \label{eq:eta1_dash}
$$
\end{theorem}

\subsection{Maximum likelihood estimation} \label{sec:mle}

\begin{theorem} \label{thm:mle}
Let $Y_1,\ldots,Y_n$ be an iid sample from the spherical Cauchy $C_d^*(\phi)$.
Suppose that ${\cal P}$ is the function (\ref{eq:ggs}).
Then the maximum likelihood estimator of $\phi$ is equal to ${\cal P}^{-1} (\hat{\theta})$, where $\hat{\theta}$ is the maximum likelihood estimator of the \sk{$d$-variate} $t$-distribution with $d$ degrees of freedom $C_d(\theta)$, given in (\ref{eq:euclid_c2}), for the sample ${\cal P}(Y_1),\ldots,{\cal P}(Y_n)$.
\end{theorem}
The proof is clear from Theorems \ref{thm:ggs} and \ref{thm:closure_two_cauchy} and omitted.
This theorem implies that, in order to estimate the parameter of the spherical Cauchy, it suffices to estimate the parameter of the \sk{$d$-variate} $t$-distribution with $d$ degrees of freedom.

\sk{Although Theorem \ref{thm:mle} is helpful for computing the maximum likelihood estimates of the parameter, there remain various properties of the maximum likelihood estimator which are not clear from this theorem.}
For example, closed-from expression of the maximum likelihood estimator for small sample size and asymptotic behaviour of the maximum likelihood estimator are not immediately obvious from Theorem \ref{thm:mle}.
The rest of this subsection is devoted to investigate properties of the maximum likelihood estimator which are not clear from Theorem \ref{thm:mle}.
The loglikelihood function is
\begin{equation}
\ell (\phi) = \sum_{j=1}^n \log f(y_j;\phi) = C + d \left\{ n \log (1-\|\phi \|^2) - \sum_{j=1}^n \log (1+\| \phi \|^2-2\phi^T y_j ) \right\},  \label{eq:loglikelihood}
\end{equation}
where $C=n \log \Gamma \{(d+1)/2 \} - n \log \{ 2 \pi^{(d+1)/2} \} $.
The first derivative of the loglikelihood function with respect to $\phi$ is
\begin{align}
\frac{\partial \ell}{\partial \phi} & =  2 d \left( \frac{n \phi }{1-\| \phi \|^2} - \sum_{j=1}^n \frac{\phi - y_j}{1+\| \phi \|^2 - 2 \phi^T y_j}  \right) = \frac{2 d}{1-\|\phi\|^2} \sum_{j=1}^n {\cal M}_{I,-\phi} (y_j), \label{eq:score}
\end{align}
\sk{where ${\cal M}_{R,\psi}$ is as in (\ref{eq:mobius_s2}).}
Therefore the estimating equation for $\phi$ has a simple form
$$
\sum_{j=1}^n {\cal M}_{I,-\phi} (y_j) = 0. \label{eq:eseq}
$$

\newpage

\begin{theorem} \label{thm:n3}
For $n \leq 3$, the maximum likelihood estimator of $\phi$, $\hat{\phi}_{ML}$, can be expressed as follows.

\begin{enumerate}
\item[(i)] For $n=1$, the maximum likelihood estimator of $\phi$ is $\hat{\phi}_{ML}=y_1$.

\item[(ii)] Suppose $n=2$.
If $y_1 \neq \pm y_2$, the contour of maximum likelihood of $\phi$ is the circle perpendicular to the unit sphere with chord $(y_1,y_2)$ in the two-dimensional plane spanned by $y_1$ and $y_2$.
When $y_1 = -y_2$, the contour of maximum likelihood of $\phi$ is the line connecting $y_1$ and $y_2$.
If $y_1=y_2$, then $\hat{\phi}_{ML}=y_1$.

\item[(iii)] Assume $n=3$ and $y_j \neq y_k$ $(j \neq k)$.
Then the maximum likelihood estimator of $\phi$ is
$$
\hat{\phi}_{ML} = {\cal P}^{-1} (\hat{\mu} + i \hat{\sigma}),
$$
where ${\cal P}$ is defined as in (\ref{eq:ggs}) and
$$
\hat{\mu} = \frac{\|{\cal P}(y_1)-{\cal P}(y_2)\|^2 {\cal P}(y_3) + \|{\cal P}(y_2)-{\cal P}(y_3)\|^2 {\cal P}(y_1) + \|{\cal P}(y_3)-{\cal P}(y_1)\|^2 {\cal P}(y_2)}{ \|{\cal P}(y_1)-{\cal P}(y_2)\|^2 + \|{\cal P}(y_2)-{\cal P}(y_3)\|^2 + \|{\cal P}(y_3)-{\cal P}(y_1)\|^2 },
$$
$$
\hat{\sigma} = \sqrt{3} \frac{\|{\cal P}(y_1)-{\cal P}(y_2)\| \|{\cal P}(y_2) - {\cal P}(y_3)\| \|{\cal P}(y_3)-{\cal P}(y_1)\|}{\| {\cal P}(y_1) - {\cal P}(y_2)\|^2 + \|{\cal P}(y_2) - {\cal P}(y_3)\|^2 + \|{\cal P}(y_3) - {\cal P}(y_1)\|^2 }.
$$
\end{enumerate}
\end{theorem}

For $d=1$ and $n=4$, \citet{mcc96} showed the maximum likelihood estimator of $\hat{\phi}_{ML}$ can be expressed in closed form.
However it does not appear clear that there are closed form expressions for the maximum likelihood estimators for $n \geq 4$ for general $d$.

\begin{lemma} \label{pro:fi}
Let $f(y)$ be the density (\ref{eq:sphere_c2}) with $\|\phi\|<1$.
Then the Fisher information matrix is
\begin{equation}
{\cal I} = - E \left\{ \frac{\partial}{\partial \phi \partial \phi^T} \log f(Y) \right\} = \frac{4}{(1-\| \phi \|^2)^2} \frac{d^2}{d+1} I. \label{eq:fisher}
\end{equation}
\end{lemma}

Thus the asymptotic variance of the maximum likelihood estimator of $\phi$ can be expressed in simple form.
\begin{theorem}
Let $Y_1,\ldots,Y_n$ be a random sample from $C^*_d(\phi)$ with $\| \phi \|<1$.
Assume $\hat{\phi}_{ML}$ is the maximum likelihood estimator of $\phi$.
Then $\sqrt{n} (\hat{\phi}_{ML} - \phi) $ tends in distribution to $N(0,{\cal I}^{-1})$ as $n \rightarrow \infty$, where ${\cal I}^{-1} = (1-\| \phi \|^2)^2 (d+1)/ (4 d^2) I$.
\end{theorem}

As seen in Theorem \ref{thm:mle}, the maximum likelihood estimates for the sample from the spherical Cauchy (\ref{eq:sphere_c2}) for general sample size can be estimated via the transformation of the spherical Cauchy into the \sk{$d$-variate} $t$ with $d$ degrees of freedom.
However it would be more efficient if the parameter estimates are obtained directly from the sample without transformation.
For $d=1$, the algorithm of \citet{ken88} is available to estimate the parameter $\phi$.
Using the fact that the Fisher information (\ref{eq:fisher}) and the score function (\ref{eq:score}) for the spherical Cauchy is expressed in simple and closed form, here we present a simple algorithm based on the Fisher scoring algorithm. \vspace{2cm}

\begin{algo} \quad \label{algo:scoring}
\begin{itemize}
\item[Step 1:\hspace{-0.4cm}] \hspace{0.3cm} Take an initial value $\phi_0$.
\item[Step 2:\hspace{-0.4cm}] \hspace{0.3cm} Compute $\phi_1,\ldots,\phi_N$ as follows until the estimate $\phi_N$ remains virtually unchanged \\
\quad \hspace{-0.05cm} from the previous estimate $\phi_{N-1}$,
$$
\phi_{t} = \phi_{t-1} + \frac{(d+1) (1-\|\phi\|^2)}{2d n} \sum_{j=1}^n \sk{ {\cal M}_{I,- \phi_{t-1}}} (y_j), \quad t=1,\ldots,N.
$$
\item[Step 3:\hspace{-0.4cm}] \hspace{0.3cm} Record $\phi_N$ as an estimate of $\phi$.
\end{itemize}
\end{algo}
The convergence of this algorithm is not proved mathematically.
However our simulation study implies that the algorithm converges fast when the method of moments estimate (\ref{eq:mme}) \sk{is used as the initial value $\phi_0$.}
In addition, for $d=1$, it seems that the parameter estimates based on Algorithm \ref{algo:scoring} numerically coincide with those based on the algorithm of \citet{ken88}.

The following tractable property holds for stationary points of the loglikelihood function.
\begin{theorem} \label{thm:unimodality}
Let $Y_1,\ldots,Y_n$ be a random sample from the spherical Cauchy $C^*(\phi)$.
Assume that $Y_j \neq Y_k$ for some $(j,k)$.
Then any stationary point of the loglikelihood function (\ref{eq:loglikelihood}) is a local maximum.
\end{theorem}


\subsection{Asymptotically efficient estimation}
Consider an estimator
\begin{equation}
\hat{\phi}_{AE} = \eta_{1,d}^{-1} (\| \overline{Y} \|) \frac{\overline{Y}}{\| \overline{Y} \|} + \frac{d+1}{2d n } \sum_{j=1}^n {\cal M}_{I,-\phi} (Y_j).  \label{eq:ae_estimator}
\end{equation}
This estimator is derived as $\hat{\phi}_{AE} = \hat{\phi}_{MM} +  (n {\cal I})^{-1} \partial \ell/ \partial \phi$, where $\hat{\phi}_{MM}$ is the method of moments estimator (\ref{eq:mme}), a consistent estimator of $\phi$, \sk{and ${\cal I}$ denotes the Fisher information matrix (\ref{eq:fisher}).}
It can be readily seen from this derivation that the estimator (\ref{eq:ae_estimator}) is an asymptotically efficient estimator of $\phi$ with asymptotic variance ${\cal I}^{-1} = (1-\| \phi \|^2)^2 (d+1)/ (4 d^2) I$.
The estimator (\ref{eq:ae_estimator}) also appears as $\phi_1$ in Algorithm \ref{algo:scoring} when \sk{the method of moments estimator (\ref{eq:mme}) is taken as the initial value $\phi_0$.}

An advantage of the estimator (\ref{eq:ae_estimator}) is that it achieves both \sk{closed-form} expression and asymptotic efficiency, whereas the method of moments estimator (\ref{eq:mme}) and the maximum likelihood estimator do not have either of them.

\section{Simulation study} \label{sec:simulation}
The method of moments estimator (\ref{eq:mme}), maximum likelihood estimator and asymptotically efficient estimator (\ref{eq:ae_estimator}) were compared in terms of the performance of finite sample sizes and the asymptotic behaviour via a Monte Carlo simulation study.
Details of the simulations given in the Supplementary Material suggest the following recommendations can be made as to the choice of the three estimators in terms of mean squared error.
If $d$ is large, then the asymptotically efficient estimator (\ref{eq:ae_estimator}) is preferred.
When $d$ is small, the asymptotically efficient estimator (\ref{eq:ae_estimator}) is preferred for dispersed data and the maximum likelihood estimator is recommended otherwise.

The calculation of the asymptotically efficient estimator (\ref{eq:ae_estimator}) is as fast as that of the method of moments estimator (\ref{eq:mme}) and is faster than that of the maximum likelihood estimator.
However the convergence of the maximum likelihood estimation based on Algorithm \ref{algo:scoring} is very fast and stable \sk{when $d$ is greater than one or $n$ is not small.}
For the circular case $d=1$, the maximum likelihood estimates estimated via Algorithm \ref{algo:scoring} numerically coincide with those estimated via the algorithm of \citet{ken88} in the sense that the sum of squared error of these two estimators is very small.

\section{Comparison with von Mises--Fisher family} \label{sec:comparison}
We compare the spherical Cauchy family with the von Mises--Fisher family which is a well-known family of distributions on the sphere.
The von Mises--Fisher family on $S^d$ has density
\begin{equation}
f(y) = \frac{\kappa^{(d-1)/2}}{(2\pi)^{(d+1)/2} I_{(d-1)/2}(\kappa)} \exp (\kappa \mu^T y ), \quad y \in S^d, \label{eq:vmf}
\end{equation}
where $\mu \in S^d$ controls the mode of the density, $\kappa \geq 0$ regulates the concentration of the distribution, and $I_{\nu}$ denotes the modified Bessel function of the first kind and order $\nu$.
The mean direction and mean resultant length of the von Mises--Fisher distribution (\ref{eq:vmf}) are $\mu$ and $A_d(\kappa)$, respectively.
See, e.g., for \citet[\S 9.3.2]{mar} for properties of von Mises--Fisher family.

First we discuss similarities and differences between the densities of the spherical Cauchy family (\ref{eq:sphere_c2}) and von Mises--Fisher family (\ref{eq:vmf}).
The densities of both families are unimodal and rotationally symmetric around their modes.
If the mean resultant lengths are small, the densities of both models have similar shapes.
In particular, when the mean resultant lengths are zero, both models reduce to the uniform distribution on the sphere.

However, when the mean resultant lengths are not small, the densities of the spherical Cauchy and von Mises--Fisher show different behaviour.
\begin{figure}
\begin{center}
\includegraphics[width=3.6cm,height=3.5cm]{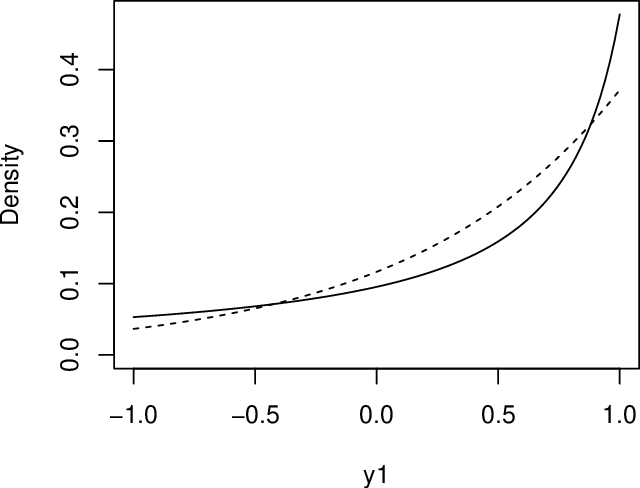} 
\includegraphics[width=3.6cm,height=3.5cm]{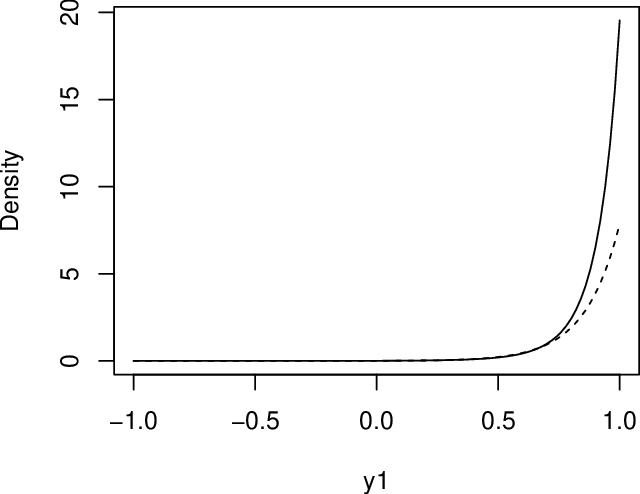} 
\includegraphics[width=3.6cm,height=3.5cm]{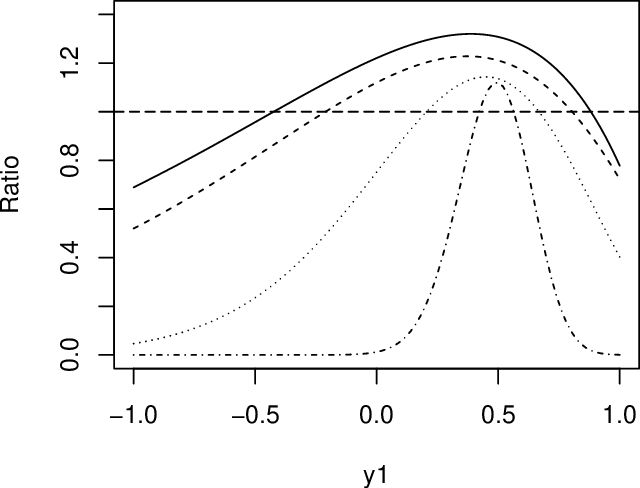} 
\includegraphics[width=3.6cm,height=3.5cm]{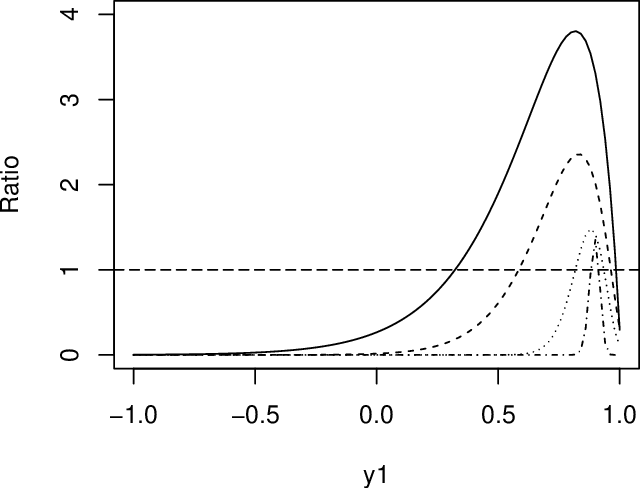}\\
\hspace{0.55cm} (a) \hspace{2.97cm} (b) \hspace{2.97cm} (c) \hspace{2.97cm} (d)
\caption{\small Density of the spherical Cauchy (\ref{eq:sphere_c2}) (solid) and that of the von Mises--Fisher (\ref{eq:vmf}) (dashed) as a function of $y_1$ for \sk{$y=(y_1, \ldots, y_{d+1})^T$,} $\phi= \eta_{1,d}^{-1}(0.5) e_1$, $\kappa \mu = A_d^{-1} (0.5) e_1$, and: (a) $d=1$ and (b) $d=10$.
\sk{The von Mises--Fisher density (\ref{eq:vmf}) divided by the spherical Cauchy density (\ref{eq:sphere_c2})} as a function of $y_1$ for $d=1$ (solid), $d=2$ (dashed), $d=10$ (dotted), and $d=100$ (dotdashed), and: (c) $\phi=\eta_{1,d}^{-1}(0.5) e_1$ and $\kappa \mu = A_d^{-1} (0.5) e_1$ and (d) $\phi=\eta_{1,d}^{-1}(0.9) e_1$ and $\kappa \mu = A_d^{-1} (0.9) e_1$.
The longdashed lines in (c) and (d) represent the horizontal lines whose intercepts are 1.
} \label{fig:comparison}
\end{center}
\end{figure}
Figure \ref{fig:comparison} displays densities and ratios of the spherical Cauchy distributions (\ref{eq:sphere_c2}) and the von Mises--Fisher distributions (\ref{eq:vmf}) for some selected values of $d$ and $\|\phi\|$.
The values of the concentration parameters are selected such that the mean resultant lengths of both models are 0.5 in Figure \ref{fig:comparison}(a)--(c) and 0.9 in Figure \ref{fig:comparison}(d).
Because of the rotational symmetry of both families, the values of the densities displayed in Figure \ref{fig:comparison} depend only on $y_1$, namely, the first component of $y$.
The figure suggests that, when the mean resultant lengths are not small, the spherical Cauchy density takes greater values than the von Mises--Fisher density around the mode and antimode and smaller values than the von Mises--Fisher density in the other area of the sphere.
The comparison between Figure \ref{fig:comparison}(a) and (b) implies that, compared with the densities with $d=2$, the densities with $d=10$ take greater values around the mode.
It seems from Figure \ref{fig:comparison}(c) and (d) that, as $d$ increases, \sk{the ratio of the two densities around the mode and that around the antimode approach zero.}
In addition, the greater the value of $d$, the smaller the range of $y_1$ in which the von Mises--Fisher density takes greater values than the spherical Cauchy density.
When the mean resultant lengths are large, the von Mises--Fisher density takes greater values than the spherical Cauchy density in a small range of $y_1$, but there is considerable difference in the values of densities in such a range. 

Next we compare other statistical aspects of the spherical Cauchy family and von Mises--Fisher family.
The von Mises--Fisher has well-developed theory of statistical inference.
Some tractable results about statistical inference for the von Mises--Fisher partly follow from the fact that, unlike the spherical Cauchy, the von Mises--Fisher is a member of the exponential family.
The maximum likelihood estimator of the parameter for the von Mises--Fisher distribution can be expressed in closed form.
On the other hand, a closed form expression for the maximum likelihood estimator has not been found apart from $n \leq 4$ for $d=1$ and $n \leq 3$ for $d \geq 2$.
As for hypothesis testing, many test statistics have been proposed in the literature for testing the location parameter and/or the concentration parameter of the von Mises--Fisher family in various settings.
Many of these test statistics are expressed in simple and closed form and their asymptotic distributions are well-studied.
Apart from the use of pivotal statistics and a direct application of likelihood ratio test, methods of hypothesis testing for the spherical Cauchy do not seem immediately clear.
Also various extensions are available for the von Mises--Fisher distribution such as the Fisher--Bingham distribution and Kent distribution \citep{ken82} and not for the spherical Cauchy distribution.

The spherical Cauchy family has the tractable property that it is closed under the M\"obius transformation on the sphere and there is a similar induced transformation on the parameter space; see Theorem \ref{thm:closure_cauchy_s}.
This result can be applied to derive tractable properties of the spherical Cauchy family such as an efficient algorithm for random variate generation, a simple form of pivotal statistics, a closed form expression for probabilities of a surface area under the spherical Cauchy density.
These properties do not hold for the von Mises--Fisher family in general.
Theorem \ref{thm:closure_cauchy_s} can also be used to simplify the computations for Fisher information matrix and maximum likelihood estimation for $n \leq 3$.
Furthermore, unlike the von Mises--Fisher family, the spherical Cauchy family is related to the $t$-family with $d$ degrees of freedom via the stereographic projection; see Theorem \ref{thm:closure_two_cauchy}.
A simple algorithm for maximum likelihood estimation and the asymptotically efficient estimator (\ref{eq:ae_estimator}) enable us to use the spherical Cauchy, which has a different shape of the density from the von Mises--Fisher in general, as a practical statistical model.
Since the M\"obius transformation and/or the wrapped Cauchy family are applied to propose statistical models for circular data including regression models and time series models, the theory of the M\"obius transformation and/or the spherical Cauchy presented in this paper can be potentially useful for the development of statistical models for spherical data.


\section*{Acknowledgment}
The first author is grateful to Department of Statistics at the University of Chicago for its hospitality during the research visit that led to this paper.
The work of the first author was supported by JSPS KAKENHI Grant Number 17K05379.

\section*{Supplementary material}
Supplementary material includes a marginal distribution of the spherical Cauchy family and its association with the real M\"obius group, moments of the marginal distribution, details of simulation study, and proofs of Lemma \ref{pro:fi} and Theorems \ref{thm:mobius_conformal}--\ref{thm:moments_sc}, \ref{thm:mme}, \ref{thm:n3} and \ref{thm:unimodality} of the article and Lemma \ref{lem:marginal} and Theorems \ref{thm:1st_moment} and \ref{thm:2nd_moment} of the supplementary material.

\newpage

\def\theequation{S\arabic{equation}}
\setcounter{equation}{0}

\def\thelemma{S\arabic{lemma}}
\setcounter{lemma}{0}

\def\thetheorem{S\arabic{theorem}}
\setcounter{theorem}{0}

\def\thesection{S\arabic{section}}
\setcounter{section}{0}

\def\thefigure{S\arabic{figure}}
\setcounter{figure}{0}

\def\thedefinition{S\arabic{definition}}
\setcounter{definition}{0}

\def\thetable{S\arabic{table}}
\setcounter{table}{0}

\title{\textbf{\Large Supplementary Material for\\\hspace{-0.14cm}``M\"obius transformation and a Cauchy family on the sphere"}\vspace{0.2cm}\\
}
\author{{\scshape Shogo Kato\(^{*,a}\) \ and Peter\ McCullagh\,\(^{b}\)} \vspace{0.5cm}\\
\textit{\(^{a}\) Institute of Statistical Mathematics, Japan} \\  
\textit{\(^{b}\) University of Chicago, USA}
}
\date{}
\maketitle

\section{A marginal distribution of a Cauchy family on $S^d$ (\S \ref{subsec:cauchy_sphere} in article)}
\subsection{A marginal distribution and real M\"obius group} \label{sec:marginal}

\begin{theorem} \label{thm:marginal}
Suppose $Y = \sk{(Y_1,\ldots,Y_{\nu+1})^T } \sim C^*_{\nu}(\phi)$, where $\phi=\sk{ (\rho,0,\ldots,0)^T }$ and $\rho \in \mathbb{R} \setminus \{-1,1\}$.
Then the marginal \sk{density} of $Y_1$ is of the form
\begin{equation}
f(y_1;\rho,\nu) = \frac{1}{B(\nu/2,1/2)} \left( \frac{|1-\rho^2|}{1+\rho^2-2 \rho y_1} \right)^{\nu} (1-y_1^2)^{(\nu-2)/2}, \quad -1 < y_1 < 1, \label{eq:marginal}
\end{equation}
where $B(\cdot,\cdot)$ denotes a beta function.
\end{theorem}
The proof is straightforward and therefore omitted.

In a similar manner as in \cite{mcc89}, if we view $\nu$ as a continuous-valued parameter with $\nu \geq 0$, then (\ref{eq:marginal}) can be considered a two-parameter family.
Clearly, $f(y_1;\rho,\nu)=f(y_1;\rho^{-1},\nu)$.
If $\rho=0$, then the distribution (\ref{eq:marginal}) reduces to the symmetric beta distribution with density
\begin{equation}
f(y_1;\nu) = \frac{(1-y_1^2)^{(\nu-2)/2}}{B(\nu/2,1/2)}, \quad -1 < y_1 <1. \label{eq:beta}
\end{equation}
It can be readily seen from equation (8.384.5) of \cite{gra} that the family (\ref{eq:marginal}) with $-1<\rho<1$ is equivalent to \citeauthor{ses}'s (\citeyear{ses}) family with the parameterization given in Example 1 of his paper.
As discussed there, if $\nu=1$, then the family (\ref{eq:marginal}) reduces to the family discussed in \cite{lei} and \cite{mcc89} whose density is given by equation (2) of the latter paper.

\begin{theorem}
Let ${\cal R}$ be the real M\"obius transformation
\begin{equation}
{\cal R}(y_1)=\frac{y_1+b}{b y_1+1}, \quad -1<y_1<1; \quad -1 < b <1. \label{eq:mobius_real}
\end{equation}
If $Y_1$ has the density (\ref{eq:marginal}), then ${\cal R}(Y_1)$ belongs to the same family with the parameter $\rho$ replaced by $(\rho+\rho')/(\rho \rho'+1)$, where $\rho'=(1-\sqrt{1-b^2})/b$.
\end{theorem}
The proof is clear from straightforward calculation and therefore omitted.
Another approach to proving this result is to remember the derivation of the model given in Theorem \ref{thm:marginal} and apply Theorem \ref{thm:closure_cauchy_s} with $R_1=R_2=I$ and $\phi_1=(\rho,0,\ldots,0)^T$ and $\phi_2=(\rho',0,\ldots,0)^T$.

This is an extension of the result given in \cite{ses} that the family (\ref{eq:marginal}) is transformed into the symmetric beta density (\ref{eq:beta}) via a special case of the M\"obius transformation (\ref{eq:mobius_real}) with $b=-2 \rho/(1+\rho^2)$.

\subsection{Moments} \label{sec:moments}

We discuss some moments of the the marginal family (\ref{eq:marginal}) which can be applied to obtain moments for the spherical Cauchy family.
Define
$$
\eta_{k,\nu} (\rho) = E(Y_1^k) = \int_{-1}^1 \frac{1}{B(\nu/2,1/2)} \left( \frac{|1-\rho^2|}{1+\rho^2-2 \rho y_1} \right)^{\nu} (1-y_1^2)^{(\nu-2)/2} dy_1, \quad -1 < \rho <1,
$$
where $Y_1$ has the density (\ref{eq:marginal}). 
As the following lemma shows, the monotonicity holds for $\eta_{k,\nu}$ for an odd integer of $k$.
\begin{lemma} \label{lem:marginal}
Suppose that $k$ is an odd integer.
Then $\eta_{k,\nu}(0) = 0$,  $\lim_{\rho \rightarrow 1} \eta_{k,\nu}(\rho) = 1$, and $ \partial \eta_{k,\nu}(\rho) / \partial \rho  > 0$ for $0 < \rho < 1$.
\end{lemma}

Next the first and second moments of the marginal family (\ref{eq:marginal}) are discussed.
\cite{ses} obtained closed-form expressions for the mean and variance of the family (\ref{eq:marginal}) with $\nu=1$ and approximated values of these statistics with general $\nu$.
Here we provide exact expressions for the moments for general $\nu$.
\begin{theorem} \label{thm:1st_moment}
The following hold for $\eta_{1,\nu} (\rho)$:
\begin{enumerate}
\item[(i)] \label{item:mu1_1}
for any $\nu \geq 0$,
\begin{align*}
\begin{split}
\eta_{1,\nu} (\rho) &= \frac{1+\rho^2}{2 \rho} \left[ 1 - \frac{(1+\rho)^2}{1+\rho^2} \, F \left\{ 1 , \frac{\nu}{2}; \nu; - \frac{4 \rho}{(1-\rho)^2} \right\} \right]  \\
 &= \frac{1+\rho^2}{2 \rho} \left[ 1 - \frac{1-\rho^2}{1+\rho^2} \, F \left\{ \frac12 , \frac{\nu-1}{2}; \frac{\nu+1}{2}; - \frac{4 \rho^2}{(1-\rho^2)^2} \right\} \right], 
\end{split}
\end{align*}
where $F$ denotes the hypergeometric series \citep[equation 9.111]{gra},

\item[(ii)] \label{item:mu1_3} for $\nu=1,\ldots,4$,
$$
\eta_{1,1} (\rho) = \rho, \quad \eta_{1,2} (\rho) = \frac{1+\rho^2}{2 \rho} \left\{ 1 - \frac{(1-\rho^2)^2}{2 \rho (1+\rho^2) } \log \left( \frac{1+\rho}{1-\rho} \right) \right\}, \quad \eta_{1,3} (\rho) = \frac{\rho (3-\rho^2)}{2},
$$
$$
\eta_{1,4} (\rho) = \frac{1+\rho^2}{2 \rho} \left\{ 1 - \frac{   3 (1-\rho^2)^2}{8 \rho^2} + \frac{3}{16 \rho^3} \frac{(1-\rho^2)^4}{1+\rho^2} \log \left( \frac{1+\rho}{1-\rho} \right)  \right\},
$$

\item[(iii)] \label{item:mu1_2}  for $\nu \geq 4$,
\sk{\begin{align*}
\eta_{1,\nu} (\rho) & = \frac{\nu-1}{(\nu-2)(\nu-3)} \Biggl[ \left\{ \nu-2 + \frac{(\nu-3) (1-\rho^2)^2}{4 \rho^2} \right\} \mu_1 (\nu-2) \\
&\ - \frac{(\nu-3)(1-\rho^2)^2}{4 \rho^2} \mu_1(\nu-4) - \frac{\nu-2}{\nu-1} \frac{1+\rho^2}{\rho} \Biggr].
\end{align*} }

\end{enumerate}
\end{theorem}

It follows from these results that, for any positive integer $\nu$, the mean of $\eta_{1,\nu} (\rho)$ can be expressed in closed form.

\begin{theorem} \label{thm:2nd_moment}
The following results hold for $\eta_{2,\nu} (\rho)$:
\begin{enumerate}
\item[(i)] \label{item:mu2_1} for $\nu \geq 0$,
\begin{align}
\label{eq:2nd_moment}
\begin{split}
\eta_{2,\nu} (\rho) = \ & \frac{(1+\rho^2)^2}{4 \rho^2} \Biggl[ 1 - 2 \frac{(1+\rho)^2}{1+\rho^2} F \left\{ 1 ,\frac{\nu}{2}; \nu; -\frac{4 \rho}{(1-\rho)^2} \right\} \\
 & + \frac{(1+\rho)^4}{(1+\rho^2)^2} F \left\{ 2 ,\frac{\nu}{2} ; \nu ; -\frac{4 \rho}{(1-\rho)^2} \right\} \Biggr],
\end{split}
\end{align}

\item[(ii)] \label{item:mu2_3}
for $\nu=1,\ldots,4$,
$$
\eta_{2,1} (\rho) = \frac{1+\rho^2}{2}, \quad \eta_{2,2} (\rho) = \frac{1+\rho^2}{4 \rho^2} \left\{ \frac{2 (1+\rho^4)}{1+\rho^2} - \frac{(1-\rho^2)^2}{\rho} \log \left( \frac{1+\rho}{1-\rho} \right) \right\},
$$
$$
\eta_{2,3} (\rho) = \frac{1+6 \rho^2 -3 \rho^4}{4},
$$
$$
\eta_{2,4} (\rho) = \frac{1+\rho^2}{16 \rho^4} \biggl\{ \frac{-2 (3 - 8 \rho^2 + 2\rho^4 -8 \rho^6 + 3 \rho^8)}{1+ \rho^2} + \frac{3 (1-\rho^2)^4}{\rho} \log \left( \frac{1+\rho}{1-\rho} \right) \biggr\},
$$

\item[(iii)] \label{item:mu2_2}
for $\nu > 4$,
\begin{align*}
\eta_{2,\nu} (\rho) &= \frac{1}{(\nu-3)(\nu-4)} \biggl[  - \frac{3(1+\rho^2)^2}{2 \rho^2} + \frac{1+\rho^2}{\rho} \bigl\{ (\nu-3)(\nu-4) \eta_{1,\nu} (\rho) \\
& - c_1 \eta_{1,\nu-2} (\rho) +  c_2 \eta_{1,\nu-4} (\rho) \bigr\} + c_1 \eta_{2,\nu-2} (\rho) - c_2 \eta_{2,\nu-4} (\rho) \biggr],
\end{align*}
where $c_1 = (\nu-1)(\nu-6) - (\nu-1)(\nu-3) (1-\rho^2)^2 /(4\rho^2)$ and $c_2 = - (\nu-1)(\nu-3) (1-\rho^2)^2/(4 \rho^2)$.
\end{enumerate}
\end{theorem}

It follows from these results and equation (9.134.3) of \cite{gra} that $\eta_{2,\nu} (\rho)$ has a \sk{closed-form} expression for any $\nu \in \mathbb{N}$. 
Thus the variance of $Y_1$ can also be expressed in closed from for any positive integer $\nu$.

Theorems \ref{thm:1st_moment} and \ref{thm:2nd_moment} can be applied to express certain moments of the spherical Cauchy family; see \S \ref{sec:mme} \sk{in the article}.


\section{Details of simulation study (\S \ref{sec:simulation} in article)}

We compare the method of moments estimator (\ref{eq:mme}), the maximum likelihood estimator and the asymptotically efficient estimator (\ref{eq:ae_estimator}) in terms of the performance of finite sample sizes and the asymptotic behaviour.
In order to compare the performance of the three estimators, the mean squared error $\textsc{MSE} = E(\| \hat{\phi} - \phi \|^2)$ is adopted, where $\hat{\phi}$ is an estimator of $\phi$ of the spherical Cauchy $C_d^*(\phi)$.
We consider the relative mean squared error defined by
$$
\textsc{RMSE}_{\textsc{X/ML}} =  \textsc{MSE}_{\textsc{X}} / \textsc{MSE}_{\textsc{ML}},
$$
where $\textsc{MSE}_{\textsc{ML}}$ denotes MSE of the maximum likelihood estimator and  $\textsc{MSE}_{\textsc{X}}$ is MSE of the method of moments estimator (\ref{eq:mme}) or the asymptotically efficient estimator (\ref{eq:ae_estimator}).

First we consider the performance of the three estimators for finite sample sizes via a Monte Carlo simulation study.
Random samples of sizes $n=10,\ 25,\ 50,\ 200$ and $1000$ were generated from the spherical Cauchy $C^*_d(\phi)$ with $\phi/\|\phi\|=e_1$, $\|\phi\|=\eta_{1,d}^{-1}(0.1),\ \eta_{1,d}^{-1}(0.3),\ \eta_{1,d}^{-1}(0.5),\ \eta_{1,d}^{-1}(0.7)$ and $\eta_{1,d}^{-1}(0.9)$ and $d=1,\ 2,\ 10,\ 50$ and $100$.
For each combination of $d$, $n$ and $\|\phi\|$, $r=2000$ random samples were generated using Corollary \ref{cor:random_variate}.
Then the three estimators were estimated for each random sample.
We used Algorithm \ref{algo:scoring} to estimate the maximum likelihood estimator \sk{and the method of moments estimator (\ref{eq:mme}) was adopted as the initial value of the algorithm}.

An estimate of MSE based on $r$ random samples is defined by $\widehat{\textsc{MSE}} = r^{-1} \sum_{j=1}^r \| \hat{\phi}_j - \phi \|^2 $, where $\hat{\phi}_j$ is an estimator of $\phi$ estimated from the $j$th random sample $(j=1,\ldots,r)$.
We then discuss an estimate of relative mean squared error defined by 
\begin{equation}
\widehat{\textsc{RMSE}}_{\textsc{X/ML}} = \widehat{\textsc{MSE}}_{\textsc{X}} \big/ \widehat{\textsc{MSE}}_{\textsc{ML}}, \label{eq:rmse_hat}
\end{equation}
where $\widehat{\textsc{MSE}}_{\textsc{ML}}$ denotes $\widehat{\textsc{MSE}}$ of the maximum likelihood estimator and $\widehat{\textsc{MSE}}_{\textsc{X}}$ is $\widehat{\textsc{MSE}}$ of the method of moments estimator (\ref{eq:mme}) or the asymptotically efficient estimator (\ref{eq:ae_estimator}).

\begin{table}
\begin{center}
\caption{\small Relative mean squared error of the method of moments estimator (\ref{eq:mme}) (MM) and that of the asymptotically efficient estimator (\ref{eq:ae_estimator}) (AE) with respect to the maximum likelihood estimator estimated from 2000 simulation samples of size $n$ from the spherical Cauchy $C_d^*(\phi)$ with $\phi/\|\phi\|=e_1$ and: (a) $d=1$, (b) $d=2$, (c) $d=10$, (d) $d=50$ and (e) $d=100$.} \vspace{0.5cm}
{\footnotesize \begin{tabular}{cccccccc}
\multicolumn{8}{c}{(a)}\vspace{0.1cm}\\
 &  & $n=10$ & $n=25$ & $n=50$ & $n=200$ & $n=1000$ & $n=\infty$ \\
\multirow{2}{*}{$\|\phi\|=\eta_{1,1}^{-1}(0.1)$} & MM & 0.914 & 0.970 & 0.987 & 1.009 & 1.005 & 1.010 \\
 & AE & 0.833 & 0.928 & 0.962 & 0.991 & 0.998 & 1.000 \\
\multirow{2}{*}{$\|\phi\|=\eta_{1,1}^{-1}(0.3)$} & MM & 0.978 & 1.055 & 1.057 & 1.103 & 1.099 & 1.099 \\
 & AE & 0.863 & 0.946 & 0.972 & 0.993 & 0.999 & 1.000 \\
\multirow{2}{*}{$\|\phi\|=\eta_{1,1}^{-1}(0.5)$} & MM & 1.153 & 1.249 & 1.303 & 1.321 & 1.326 & 1.333 \\
 & AE & 0.927 & 0.984 & 0.994 & 0.995 & 1.000 & 1.000 \\
\multirow{2}{*}{$\|\phi\|=\eta_{1,1}^{-1}(0.7)$} & MM & 1.603 & 1.849 & 1.832 & 1.902 & 1.931 & 1.961 \\
 & AE & 1.108 & 1.074 & 1.041 & 1.013 & 1.004 & 1.000 \\
\multirow{2}{*}{$\|\phi\|=\eta_{1,1}^{-1}(0.9)$} & MM & 3.861 & 4.657 & 4.850 & 4.917 & 5.286 & 5.263 \\
 & AE & 2.128 & 1.902 & 1.545 & 1.151 & 1.036 & 1.000 \\
 \\
\multicolumn{8}{c}{(b)}\vspace{0.1cm}\\
 & & $n=10$ & $n=25$ & $n=50$ & $n=200$ & $n=1000$ & $n=\infty$ \\
\multirow{2}{*}{$\|\phi\|=\eta_{1,2}^{-1}(0.1)$} & MM & 0.965 & 0.985 & 1.002 & 1.000 & 1.003 & 1.005 \\
 & AE & 0.903 & 0.961 & 0.980 & 0.995 & 0.999 & 1.000 \\
\multirow{2}{*}{$\|\phi\|=\eta_{1,2}^{-1}(0.3)$} & MM & 0.996 & 1.027 & 1.043 & 1.046 & 1.043 & 1.048 \\
 & AE & 0.903 & 0.962 & 0.983 & 0.996 & 0.999 & 1.000 \\
\multirow{2}{*}{$\|\phi\|=\eta_{1,2}^{-1}(0.5)$} & MM & 1.096 & 1.123 & 1.141 & 1.158 & 1.159 & 1.153 \\
 & AE & 0.923 & 0.970 & 0.987 & 0.997 & 1.000 & 1.000 \\
\multirow{2}{*}{$\|\phi\|=\eta_{1,2}^{-1}(0.7)$} & MM & 1.288 & 1.336 & 1.377 & 1.369 & 1.415 & 1.392 \\
 & AE & 0.961 & 0.986 & 0.997 & 0.999 & 1.000 & 1.000 \\
\multirow{2}{*}{$\|\phi\|=\eta_{1,2}^{-1}(0.9)$} & MM & 1.896 & 2.165 & 2.159 & 2.222 & 2.275 & 2.234 \\
 & AE & 1.071 & 1.043 & 1.031 & 1.011 & 1.002 & 1.000 \\
 \\
\multicolumn{8}{c}{(c)}\vspace{0.1cm}\\
 & & $n=10$ & $n=25$ & $n=50$ & $n=200$ & $n=1000$ & $n=\infty$ \\
\multirow{2}{*}{$\|\phi\|=\eta_{1,10}^{-1}(0.1)$} & MM & 0.993 & 0.996 & 0.999 & 1.001 & 1.002 & 1.001 \\
 & AE & 0.977 & 0.991 & 0.996 & 0.999 & 1.000 & 1.000 \\
\multirow{2}{*}{$\|\phi\|=\eta_{1,10}^{-1}(0.3)$} & MM & 0.998 & 1.002 & 1.005 & 1.010 & 1.009 & 1.009 \\
 & AE & 0.977 & 0.991 & 0.996 & 0.999 & 1.000 & 1.000 \\
\multirow{2}{*}{$\|\phi\|=\eta_{1,10}^{-1}(0.5)$} & MM & 1.013 & 1.024 & 1.023 & 1.024 & 1.031 & 1.027 \\
 & AE & 0.978 & 0.992 & 0.996 & 0.999 & 1.000 & 1.000 \\
\multirow{2}{*}{$\|\phi\|=\eta_{1,10}^{-1}(0.7)$} & MM & 1.037 & 1.050 & 1.062 & 1.060 & 1.056 & 1.058 \\
 & AE & 0.980 & 0.993 & 0.996 & 0.999 & 1.000 & 1.000 \\
\multirow{2}{*}{$\|\phi\|=\eta_{1,10}^{-1}(0.9)$} & MM & 1.084 & 1.106 & 1.111 & 1.106 & 1.110 & 1.111 \\
 & AE & 0.984 & 0.994 & 0.997 & 0.999 & 1.000 & 1.000 \\
  \\
\multicolumn{8}{c}{(d)}\vspace{0.1cm}\\
 & & $n=10$ & $n=25$ & $n=50$ & $n=200$ & $n=1000$ & $n=\infty$ \\
\multirow{2}{*}{$\|\phi\|=\eta_{1,50}^{-1}(0.1)$} & MM & 0.998 & 1.000 & 1.000 & 1.000 & 1.000 & 1.000 \\
 & AE & 0.995 & 0.998 & 0.999 & 1.000 & 1.000 & 1.000 \\
\multirow{2}{*}{$\|\phi\|=\eta_{1,50}^{-1}(0.3)$} & MM & 0.999 & 1.001 & 1.001 & 1.002 & 1.001 & 1.002 \\
 & AE & 0.995 & 0.998 & 0.999 & 1.000 & 1.000 & 1.000 \\
\multirow{2}{*}{$\|\phi\|=\eta_{1,50}^{-1}(0.5)$} & MM & 1.003 & 1.004 & 1.005 & 1.005 & 1.004 & 1.005 \\
 & AE & 0.995 & 0.998 & 0.999 & 1.000 & 1.000 & 1.000 \\
\multirow{2}{*}{$\|\phi\|=\eta_{1,50}^{-1}(0.7)$} & MM & 1.007 & 1.008 & 1.009 & 1.010 & 1.011 & 1.010 \\
 & AE & 0.995 & 0.998 & 0.999 & 1.000 & 1.000 & 1.000 \\
\multirow{2}{*}{$\|\phi\|=\eta_{1,50}^{-1}(0.9)$} & MM & 1.014 & 1.016 & 1.017 & 1.016 & 1.019 & 1.017 \\
 & AE & 0.996 & 0.999 & 0.999 & 1.000 & 1.000 & 1.000
\end{tabular}
}
\label{tab:estimators}
\end{center}
\end{table}

\begin{table}
\begin{center}
{\footnotesize 
\begin{tabular}{cccccccc}
\multicolumn{8}{c}{(e)}\vspace{0.1cm}\\
 & & $n=10$ & $n=25$ & $n=50$ & $n=200$ & $n=1000$ & $n=\infty$ \\
\multirow{2}{*}{$\|\phi\|=\eta_{1,100}^{-1}(0.1)$} & MM & 0.998 & 1.000 & 1.000 & 1.000 & 1.000 & 1.000 \\
 & AE & 0.998 & 0.999 & 1.000 & 1.000 & 1.000 & 1.000 \\
\multirow{2}{*}{$\|\phi\|=\eta_{1,100}^{-1}(0.3)$} & MM & 1.000 & 1.000 & 1.001 & 1.001 & 1.001 & 1.001 \\
 & AE & 0.998 & 0.999 & 1.000 & 1.000 & 1.000 & 1.000 \\
\multirow{2}{*}{$\|\phi\|=\eta_{1,100}^{-1}(0.5)$} & MM & 1.001 & 1.003 & 1.003 & 1.002 & 1.002 & 1.003 \\
 & AE & 0.998 & 0.999 & 1.000 & 1.000 & 1.000 & 1.000 \\
\multirow{2}{*}{$\|\phi\|=\eta_{1,100}^{-1}(0.7)$} & MM & 1.004 & 1.005 & 1.005 & 1.005 & 1.004 & 1.005 \\
 & AE & 0.998 & 0.999 & 1.000 & 1.000 & 1.000 & 1.000 \\
\multirow{2}{*}{$\|\phi\|=\eta_{1,100}^{-1}(0.9)$} & MM & 1.007 & 1.008 & 1.008 & 1.008 & 1.008 & 1.008 \\
 & AE & 0.998 & 0.999 & 1.000 & 1.000 & 1.000 & 1.000
\end{tabular}
}
\end{center}
\end{table}

Table \ref{tab:estimators} shows estimates of relative mean squared error (\ref{eq:rmse_hat}) for some selected combinations of $d$, $n$ and $\|\phi\|$.
The values of $\|\phi\|$ are defined such that the mean resultant lengths of the underlying distributions are $0.1$, $0.3$, $0.5$, $0.7$ and \sk{$0.9$.}
The values of the relative mean squared error (\ref{eq:rmse_hat}) for $n=\infty$ given in the table are those of the asymptotic relative mean squared error, namely, $\lim_{n \rightarrow \infty} \textsc{RMSE}_{\textsc{X/ML}} $, which can be calculated using Theorem \ref{thm:mme} and Lemma \ref{pro:fi}.

The conclusions deduced from this table are the following.
For high dimensional cases, that is, $d \geq 10$, the asymptotically efficient estimator (\ref{eq:ae_estimator}) slightly outperforms the method of moments estimator (\ref{eq:mme}) and the maximum likelihood estimator in terms of the mean squared error.
For low dimensional cases, that is, $d=1$ or $2$, the asymptotically efficient estimator (\ref{eq:ae_estimator}) outperforms the other two estimators for small values of $\|\phi\|$ and the maximum likelihood estimator is preferable otherwise.
The method of moments estimator (\ref{eq:mme}) shows worse performance than the asymptotically efficient estimator (\ref{eq:ae_estimator}) in all the cases, especially, those of small $d$ and large $\|\phi\|$.

In a little more detail, both the asymptotically efficient estimator (\ref{eq:ae_estimator}) and the method of moments estimator (\ref{eq:mme}) outperform the maximum likelihood estimator in terms of mean squared error in the cases of small $n$ and $\|\phi\|$.
In those cases, the smaller the value of $d$, the better the performance of the asymptotically efficient estimator (\ref{eq:ae_estimator}) over the other two estimators.

For fixed values of $d$ and $\|\phi\|$, as $n$ increases, the value of $\widehat{\textsc{RMSE}}_{\textsc{X/ML}}$ for the method of moments estimator (\ref{eq:mme}) increases and the value of $\widehat{\textsc{RMSE}}_{\textsc{X/ML}}$ for the asymptotically efficient estimator (\ref{eq:ae_estimator}) approaches one.
For fixed $d$ and $n$, the greater the value of $\|\phi\|$, the greater the value of $\widehat{\textsc{RMSE}}_{\textsc{X/ML}}$ for both the method of moments estimator (\ref{eq:mme}) and the asymptotically efficient estimator (\ref{eq:ae_estimator}).
The values of $\widehat{\textsc{RMSE}}_{\textsc{X/ML}}$ for the method of moments estimator (\ref{eq:mme}) are greater than those for the asymptotically efficient estimator (\ref{eq:ae_estimator}) \sk{in all the cases}.

These trends of the performance among the three estimators given in the last paragraph are particularly clear for small $d$.
As $d$ increases, the difference among the estimators in terms of mean squared error decreases.
In particular, when $d=100$, there are not considerable differences among the three estimators.

\begin{figure}
\begin{center}
\begin{tabular}{cc}
\includegraphics[width=6cm,height=5cm]{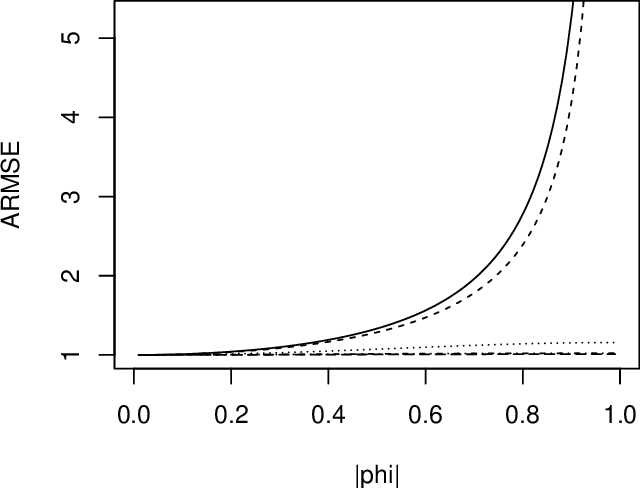} &
\includegraphics[width=6cm,height=5cm]{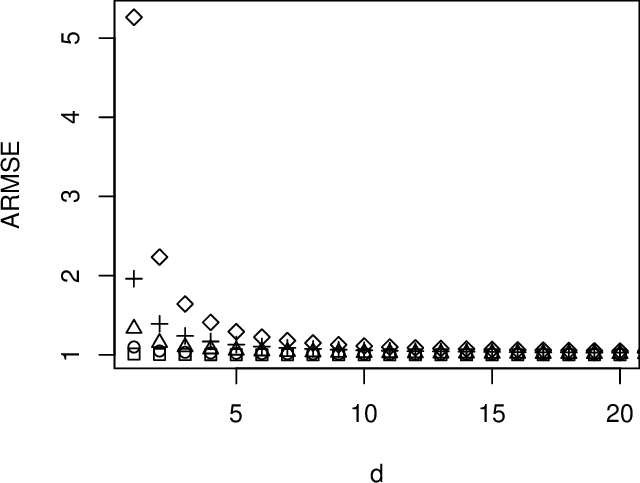} \\
\hspace{1.1cm} {\small (a)} & \hspace{1.05cm} {\small (b)}
\end{tabular}
\caption{\small Asymptotic relative mean squared error of the method of moments estimator (\ref{eq:mme}) with respect to the maximum likelihood estimator for the spherical Cauchy $C_d^*(\phi)$ as a function of: (a) $\| \phi \|$ for $d=1$ (solid), $d=2$ (dashed), $d=10$ (dotted), $d=50$ (dotdahsed), and $d=100$ (longdashed) and (b) $d$ for $\| \phi \|=\eta_{1,d}^{-1}(0.1) (\Box)$, $\| \phi \|=\eta_{1,d}^{-1}(0.3) (\bigcirc )$, $\| \phi \|=\eta_{1,d}^{-1}(0.5) (\bigtriangleup)$, $\| \phi \|=\eta_{1,d}^{-1}(0.7) (+)$, and $\| \phi \|=\eta_{1,d}^{-1}(0.9) (\Diamond)$.} \label{fig:armse}
\end{center}
\end{figure}

Table \ref{tab:estimators} suggests that the values of the relative mean squared error of the method of moments estimator (\ref{eq:mme}) with respect to the maximum likelihood estimator increase as $n$ increases.
Here we discuss more details on the limits of these values as $n \rightarrow \infty$.
Figure \ref{fig:armse} displays the asymptotic relative mean squared error of the method of moments estimator (\ref{eq:mme}) with respect to the maximum likelihood estimator as a function of $d$ or $n$.
Its panel (a) implies that, when $\|\phi\|$ is small, the asymptotic relative mean squared error is close to one for any $d$.
This panel also suggests that the asymptotic relative mean squared error is monotonically increasing with respect to $\|\phi\|$.
In particular, when $d$ is small and $\|\phi\|$ is large, the asymptotic relative mean squared error is very large.
As $d$ increases, The asymptotic relative mean squared error approaches one for any $\|\phi\|$.
Figure \ref{fig:armse}(b) implies that, when the mean resultant length is small, the asymptotic mean squared error of the method of moments estimator (\ref{eq:mme}) is close to that of the maximum likelihood estimator.
If the mean resultant length is large, the asymptotic relative mean squared error is large for small $d$ and close to one for large $d$.
The asymptotic relative mean squared error monotonically decreases as $d$ increases.
\sk{Because of asymptotically efficiency of the asymptotically efficient estimator (\ref{eq:ae_estimator}), the same discussion can be given to the relative mean squared error of the method of moments estimator (\ref{eq:mme}) with respect to the asymptotically efficient estimator.}

Given these observations, the following conclusions can be made as to the choice of the three estimators of the parameter of the spherical Cauchy in terms of mean squared error.
If the dimension of the data is large, then the asymptotically efficient estimator (\ref{eq:ae_estimator}) is preferred.
This estimator outperforms both the maximum likelihood estimator and method of moments estimator (\ref{eq:mme}) in terms of mean squared error for large $d$.
When the dimension of the data is small, the asymptotically efficient estimator (\ref{eq:ae_estimator}) is preferred for dispersed data and the maximum likelihood estimator is recommended otherwise.

The calculation of the asymptotically efficient estimator (\ref{eq:ae_estimator}) is as efficient as that of the method of moments estimator (\ref{eq:mme}) and is more efficient than that of the maximum likelihood estimator.
However, our simulation study suggests that the converge of the maximum likelihood estimation based on Algorithm \ref{algo:scoring} is very fast when $n$ is not very small and $d$ is greater than one.
Actually, our computation for producing Table \ref{tab:estimators} implies that Algorithm \ref{algo:scoring} converges in almost all the combinations of $(d,n,\|\phi\|)$ \sk{when the method of moments estimator (\ref{eq:mme}) is adopted as the initial value}.
To be more precise, using the stopping rule $\|\phi_t - \phi_{t-1}\| < 1 \times 10^{-7} $ and $t \leq 100$, Algorithm \ref{algo:scoring} failed to converge only twice for $(d,n,\|\phi\|)=(1,10,\eta_{1,1}^{-1}(0.1))$ and once for $(d,n,\|\phi\|)=(1,10,\eta_{1,1}^{-1}(0.7))$ and $(d,n,\|\phi\|)=(1,10,\eta_{1,1}^{-1}(0.9))$ among 2000 simulation samples for each combination of $(d,n,\|\phi\|)$.
When the stopping rule is relaxed to be $\|\phi_t - \phi_{t-1}\| < 1 \times 10^{-5} $ and $t \leq 100$, then Algorithm \ref{algo:scoring} converged in all the cases.
Also, when $d=1$, our simulation study implies that the maximum likelihood estimates estimated via Algorithm \ref{algo:scoring} numerically coincide with those estimated via the algorithm of \citet{ken88} in the sense that the sum of squared error of these two estimators is very small.

\section{Proofs}
\subsection{Proof of Theorem \ref{thm:mobius_conformal}}

\begin{proof}
\begin{enumerate}
\item[(i)]
Consider the function
\begin{equation}
\check{\cal M}(x) = A \left( \gamma \, \frac{x+a}{\|x+a\|^{\varepsilon}} + b \right), \quad x \in \mathbb{R}^{d+1} \setminus \{-a\}, \label{eq:mobius_r}
\end{equation}
where $a,b \in \mathbb{R}^{d+1}$, \sk{$\gamma \in \mathbb{R} \setminus \{ 0 \}$}, $A$ is a $(d+1) \times (d+1)$ orthogonal matrix, and $\varepsilon$ is either 0 or 2.
If $x \in \{-a,\infty\}$, we define $\check{\cal M}(-a)=Ab$ and $\check{\cal M}(\infty)=\infty$ for $\varepsilon=0$ and $\check{\cal M}(-a)=\infty$ and $\check{\cal M}(\infty) = b$ for $\varepsilon=2$.
It is known that the transformation (\ref{eq:mobius_r}) is a bijective conformal map which maps $\overline{\mathbb{R}}^{d+1}$ onto itself \citep[see][\S 2]{iwa}.
Since the transformation (\ref{eq:mobius_s2}), which has the alternative expression (\ref{eq:mobius_s3}), is a special case of (\ref{eq:mobius_r}), it follows that (\ref{eq:mobius_s2}) is a bijective conformal map which maps $\overline{\mathbb{R}}^{d+1}$ onto itself.

\item[(ii)] This is clear from the fact the transformation (\ref{eq:mobius_s2}) reduces to (\ref{eq:mobius_s}) if $x \in S^d$.

\item[(iii)] It holds that
\begin{align}
\| {\cal M}_{R,\sk{\psi}}(x) \|^2 &= \left\| R \, T_{\sk{\psi}} \left\{ \frac{1- \| \tilde{\sk{\psi}} \|^2}{\| x + \tilde{\sk{\psi}} \|^2} ( x+ \tilde{\sk{\psi}}) + \tilde{\sk{\psi}} \right\} \right\|^2  = \frac{1+2 x^T \tilde{\sk{\psi}} + \|x\|^2 \|\tilde{\sk{\psi}}\|^2 }{ \| x+\tilde{\sk{\psi}} \|^2 }. \label{eq:norm_g}
\end{align}
The difference between the numerator and denominator of the last expression of (\ref{eq:norm_g}) is
$$
(1+2 x^T \tilde{\sk{\psi}} + \|x\|^2 \|\tilde{\sk{\psi}}\|^2) - \| x+\tilde{\sk{\psi}} \|^2 = (1-\|x\|^2) (1-\|\tilde{\sk{\psi} }\|^2).
$$
This implies that, if $\|\sk{\psi}\|<1$, then $ {\cal M}_{R,\sk{\psi}} (x) \in D^{d+1}$ for $x \in D^{d+1}$ and $ {\cal M}_{R,\sk{\psi}} (x) \in \overline{\mathbb{R}}^{d+1} \setminus \overline{D}^{d+1} $ for $x \in \overline{\mathbb{R}}^{d+1} \setminus \overline{D}^{d+1}$.
It follows from this fact and the bijectivity of $ {\cal M}_{R,\sk{\psi}}$ given in (i) that, for any $y \in D^{d+1}$, there exists $x \in D^{d+1}$ such that $ y={\cal M}_{R,\sk{\psi}} (x) $.
Similarly, for any $y \in \overline{\mathbb{R}}^{d+1} \setminus \overline{D}^{d+1} $, there exists $x \in \overline{\mathbb{R}}^{d+1} \setminus \overline{D}^{d+1} $ such that $ y={\cal M}_{R,\sk{\psi}} (x) $.

\item[(iv)]
This result can be proved in a similar manner as in (iii).
\end{enumerate}
\end{proof}

\subsection{Proof of Theorem \ref{thm:group}}

\begin{proof}
For convenience write ${\cal M}_{R_2,\sk{\psi}_2} \circ {\cal M}_{R_1,\sk{\psi}_1} (x) = {\cal M}_{R_2,\sk{\psi}_2} \{ {\cal M}_{R_1,\sk{\psi}_1} (x)    \}$ for ${\cal M}_{R_1,\sk{\psi}_1},{\cal M}_{R_2,\sk{\psi}_2} \in {\cal G}$.
Then Lemma \ref{lem:closure_mobius_s} implies that, if ${\cal M}_{R_1,\sk{\psi}_1}, {\cal M}_{R_2,\sk{\psi}_2} \in {\cal G}$, then ${\cal M}_{R_2,\sk{\psi}_2} \circ {\cal M}_{R_1,\sk{\psi}_1} \in {\cal G}$.
Also, it is clear that $({\cal M}_{R_3,\sk{\psi}_3} \circ {\cal M}_{R_2,\sk{\psi}_2}) \circ {\cal M}_{R_1,\sk{\psi}_1} =  {\cal M}_{R_3,\sk{\psi}_3} \circ ( {\cal M}_{R_2,\sk{\psi}_2} \circ {\cal M}_{R_1,\sk{\psi}_1} )$.
It follows from Lemma \ref{lem:closure_mobius_s} that, for any ${\cal M}_{R,\sk{\psi}} \in {\cal G}$, ${\cal M}_{R,\sk{\psi}} \circ {\cal M}_{I,0} = {\cal M}_{I,0} \circ {\cal M}_{R,\sk{\psi}} = {\cal M}_{R,\sk{\psi}}$, implying that ${\cal M}_{I,0} \in {\cal G}$ is the identity element of ${\cal G}$.
Finally, the existence of the inverse element of ${\cal G}$ can be seen from the fact that, for any ${\cal M}_{R,\sk{\psi}} \in {\cal G}$, it holds that ${\cal M}_{R^T, - R \sk{\psi}}  \in {\cal G}$ and ${\cal M}_{R,\sk{\psi}} \circ {\cal M}_{R^T, - R \sk{\psi}} = {\cal M}_{R^T, - R \sk{\psi}} \circ {\cal M}_{R,\sk{\psi}} = {\cal M}_{I,0}$. 
\end{proof}

\subsection{Proof of Theorem \ref{thm:closure_cauchy_s}}

\begin{proof}
\sk{Let $U=(U_1, \ldots, U_{d+1})^T$ follow the uniform distribution on the sphere with density $f(u)= \Gamma \{(d+1)/2\}/\{2 \pi^{(d+1)/2}\}$.
Suppose $U_1 = \cos \tilde{\Theta}_1$, $ U_j = (\prod_{k=1}^{j-1} \sin \tilde{\Theta}_k ) \cos \tilde{\Theta}_j$ and $U_{d+1} = \prod_{k=1}^d \sin \tilde{\Theta}_k$ $(2 \leq j \leq d)$, where $\tilde{\Theta}_1, \ldots, \tilde{\Theta}_{d-1}$ take values in $[0,\pi]$ and $\tilde{\Theta}_d$ takes values in $[0,2\pi)$.
Then the density of $(\tilde{\Theta}_1,\ldots,\tilde{\Theta}_d)$ is
$$
f(\tilde{\theta}_1,\ldots,\tilde{\theta}_d) = \frac{\Gamma \{(d+1)/2\} }{ 2 \pi^{(d+1)/2} } \sin^{d-1} \tilde{\theta}_1 \sin^{d-2} \tilde{\theta}_2 \cdots \sin \tilde{\theta}_{d-1}.
$$
Next, let $V= (V_1,\ldots,V_{d+1}) = {\cal M}_{I,\rho e_1}(U)$, where $0 \leq \rho <1$.
Assume $V_1 = \cos \Theta_1$, $ V_j = (\prod_{k=1}^{j-1} \sin \Theta_k ) \cos \Theta_j$ and $V_{d+1} = \prod_{k=1}^d \sin \Theta_k$ $(2 \leq j \leq d)$, where $\Theta_1, \ldots, \Theta_{d-1}$ take values in $[0,\pi]$ and $\Theta_d$ takes values in $[0,2\pi)$.
Then it holds that $\cos \Theta_1 = (1-\rho^2)/(1+\rho^2-2 \rho \cos \tilde{\Theta}_1) - \rho$ and $\Theta_j = \tilde{\Theta}_j$ $(2 \leq j \leq d)$.
It follows that
$$
f(\theta_1,\ldots,\theta_d) =  \frac{\Gamma \{(d+1)/2\} }{ 2 \pi^{(d+1)/2} } \left( \frac{1-\rho^2}{1+\rho^2-2\rho \cos \theta_1} \right)^d \sin^{d-1} \theta_1 \sin^{d-2} \theta_2 \cdots \sin \theta_{d-1},
$$
implying $V \sim C^*(\rho e_1)$.
Since the uniform distribution on the sphere is invariant under rotation, we have ${\cal M}_{R,\phi}(U) \sim C^*_d (R \phi)$ for any $R$ and $\phi$.}

Also, it follows from Lemma \ref{lem:closure_mobius_s} that ${\cal M}_{R,\psi} \{ {\cal M}_{I,\phi} (U) \} = {\cal M}_{\check{R},\check{\phi}} (U)$, where $\check{R}=R T_{\psi} T_{\beta} T_{\phi} T_{\check{\phi}}$, $\check{\phi} = T_{\phi} T_{\beta} T_{\psi} \sk{ {\cal M}_{I,\psi} } (\phi)$ and $\beta=\tilde{\phi} + \tilde{\psi}$.
Therefore ${\cal M}_{I,\phi} (U) (\equiv Y)$ has the spherical Cauchy distribution $C^*_d (\phi)$ and its transformation ${\cal M}_{R,\psi} ( Y )$ also follows the spherical Cauchy $C^*_d (\omega)$, where
\begin{align*}
\omega &= \check{R} \check{\phi} =  R T_{\psi} T_{\beta} T_{\phi} \cdot  T_{\phi} T_{\beta} T_{\psi} \sk{ {\cal M}_{I,\psi} (\phi) } =  {\cal M}_{R,\psi} (\phi).
\end{align*}
\end{proof}

\subsection{Proof of Theorem \ref{thm:ggs}}

\begin{proof}
We first prove (i).
Consider the function
\begin{equation}
{\cal P}^\star(x) = \left( I - 2 e_{d+1} e^T_{d+1} \right) \left\{ 2 \, \frac{x-e_{d+1}}{\|x-e_{d+1}\|^2} + e_{d+1} \right\}, \quad x \in \mathbb{R}^{d+1} \setminus \{e_{d+1}\}. \label{eq:ggs2}
\end{equation}
Also, assume that ${\cal P}^\star(\infty) = -e_{d+1}$ and ${\cal P}^\star(e_{d+1}) = \infty$.
Clearly the function (\ref{eq:ggs2}) is a M\"obius transformation on $\overline{\mathbb{R}}^{d+1}$ which maps $\overline{\mathbb{R}}^{d+1}$ onto itself; see \citet[\S 2]{iwa}.
It can be seen that ${\cal P}(x)$ is equal to ${\cal P}^\star(x)$ if the imaginary part of ${\cal P}(x)$ is identified as the $(d+1)$-th component of ${\cal P}^\star(x)$.
It follows from the \sk{general} property of the M\"obius transformation that the function (\ref{eq:ggs}) is a bijective function which maps $(\mathbb{R}^d + i \mathbb{R}) \cup \{\infty\}$ onto $\overline{\mathbb{R}}^{d+1}$.
The other properties (ii) and (iii) are clear from the definition of \sk{${\cal P}$}.
\end{proof}

\subsection{Proof of Theorem \ref{thm:closure_two_cauchy}}

\begin{proof}
For convenience, write $\tilde{Y} = (\tilde{Y}_1, \ldots, \tilde{Y}_{d+1}) = (Y_{d+1},\ldots,Y_1)$.
Let $\tilde{Y}_1 = \cos \Theta_1$, $ \tilde{Y}_j = (\prod_{k=1}^{j-1} \sin \Theta_k ) \cos \Theta_j$ \sk{and} $\tilde{Y}_{d+1} = \prod_{k=1}^d \sin \Theta_k$ $(2 \leq j \leq d)$, where $\Theta_1, \ldots, \Theta_{d-1}$ take values in $[0,\pi]$ and $\Theta_d$ takes values in $[0,2\pi)$.
Then the density of $(\Theta_1,\ldots,\Theta_d)$ is of the form
$$
f(\theta_1,\ldots,\theta_d) = \frac{\Gamma \{(d+1)/2\}}{2 \pi^{(d+1)/2}} \left( \frac{ \left| 1- \|\phi \|^2 \right| }{\| y(\theta) - \phi \|^2} \right)^d \sin^{d-1} \theta_1 \sin^{d-2} \theta_2 \cdots \sin \theta_{d-1},
$$
where $y(\theta) = (\tilde{y}_{d+1},\ldots,\tilde{y}_1)^T$ is a function of $\theta_1,\ldots,\theta_d$.

Using $(\theta_1,\ldots,\theta_d)$, the function (\ref{eq:ggs}) defined on the sphere or stereographic projection (\ref{eq:gis}) can be expressed as
$$
\tilde{\cal P} (\theta_1,\ldots,\theta_d) = \frac{\sin \theta_1}{1-\cos \theta_1} \left\{ \prod_{k=2}^{d-1} \sin \theta_k, \left( \prod_{k=2}^{d-2} \sin \theta_k \right) \cos \theta_{d-1}, \ldots , \sin \theta_2 \cos \theta_3, \cos \theta_2 \right\}.
$$
Write $z = (z_1,\ldots,z_d) = \tilde{\cal P} (\theta_1,\ldots,\theta_d) $.
Then the density of $z$ is of the form
\begin{align*}
f(z) & = \frac{\Gamma \{(d+1)/2\}}{2 \pi^{(d+1)/2}} \left( \frac{ \left| 1- \|\phi \|^2 \right| }{\| y(z) - \phi \|^2} \right)^d \left( \frac{2}{1+\|z\|^2} \right)^d   \\
 & = \frac{2^{d-1} \Gamma \{(d+1)/2\}}{\pi^{(d+1)/2}} \left( \frac{ |\sigma_{\phi} | }{ \sigma_{\phi}^2 + \| z - \mu_{\phi} \|^2} \right)^d,
\end{align*}
where $y(z)=\{ 2/(1+\|z\|^2) \} ( z_1,\ldots,z_d,(\|z\|^2-1)/2 )$ and $\mu_{\phi} + i \sigma_{\phi}={\cal P}(\phi)$.

Since ${\cal P}$ is a bijective mapping, it is easy to see that the multivariate $t$-distribution $C_d(\theta)$ is transformed into the spherical Cauchy $C_d^*\{ {\cal P}^{-1}(\theta)\}$ via the inverse stereographic projection ${\cal P}^{-1}$.
\end{proof}

\subsection{Proof of Theorem \ref{thm:moments_sc}}

\begin{proof}
We first prove the case $\phi \neq 0$.
Let $\tilde{Y}=(\tilde{Y}_1,\ldots,\tilde{Y}_{d+1})^T$ follow the spherical Cauchy $C_d^*(\tilde{\phi})$, where $\tilde{\phi}=(\|\phi\|,0,\ldots,0)^T$.
Then it follows from Theorem \ref{thm:1st_moment} that $E(\tilde{Y}_1)= \eta_{1,d} (\|\phi\|)$.
Also, the symmetry of the marginal density of $\tilde{Y}_j$ $(2 \leq j \leq d+1)$ implies $E(\tilde{Y}_j) = 0$.
Then we obtain $E(Y)$ by transforming $\tilde{Y}$ via $Y = R \tilde{Y}$, where $R$ is a rotation matrix whose first column is $\phi / \|\phi\|$.

As for $E(Y Y^T)$, Theorem \ref{thm:2nd_moment} implies $E(\tilde{Y}_1^2) = \eta_{2,d} (\|\phi\|)$.
In order to calculate the moment $E(\tilde{Y}_2^2)$ \sk{for $d \geq 3$}, we first transform $\tilde{Y}$ into polar-coordinate form such that the first and second elements of $\tilde{Y}$ are $\tilde{Y}_1 = \cos \Theta_1$ and $\tilde{Y}_2 = \sin \Theta_1 \cos \Theta_2$, respectively, where $\Theta_1$ and $\Theta_2$ take values in $[0,\pi]$.
Then $E(\tilde{Y}_2^2)$ is of the form
$$
E(\tilde{Y}_2^2) = \frac{d-1}{2\pi} (1- \|\phi\|^2)^d \int_0^{\pi} \frac{\sin^{d+1} \theta_2}{(1+\|\phi\|^2-2\|\phi\| \cos \theta_1)^d} d \theta_1 \int_0^{\pi} \cos^2 \theta_2 \sin^{d-2} \theta_2 d\theta_2.
$$
Using Theorem \ref{thm:2nd_moment}, we have $E(\tilde{Y}_2^2) = d^{-1} \{ 1- \eta_{2,d} (\|\phi\|) \}$.
\sk{The moment $E(\tilde{Y}_2^2)$ for $d =1$ and $d=2$ can be calculated in a similar manner by transforming $(\tilde{Y}_1,\tilde{Y}_2)$ into polar-coordinate form.}
As for \sk{$E(\tilde{Y}_j \tilde{Y}_k)$} $(j \neq k)$, the symmetry of the marginal distribution \sk{$(\tilde{Y}_j, \tilde{Y}_k)$} implies \sk{$E(\tilde{Y}_j \tilde{Y}_k)=0$.}
Then $E(\tilde{Y} \tilde{Y}^T) = \mbox{diag} [ \eta_{2,d} (\|\phi\|), d^{-1} \{ 1- \eta_{2,d} (\|\phi\|) \}, \ldots, d^{-1} \{1-\eta_{2,d} (\|\phi\|) \} ]$.
Transforming $\tilde{Y}$ via $Y = R \tilde{Y}$, where $R$ is a rotation matrix whose first column is $\phi / \|\phi\|$, we obtain $E(YY^T)$.

If $\phi = 0$, then $Y$ has the uniform distribution on the sphere.
In this case it is known that $E(Y)=0$ and $E(Y Y^T) =(d+1)^{-1} I$; see \citet[\S 9.6.1]{mar}.
\end{proof}

\subsection{Proof of Theorem \ref{thm:mme}}

\begin{proof}
It follows from Theorem \ref{thm:moments_sc} that $\sqrt{n}(\overline{Y} - \eta_{1,d} (\phi) \phi/\|\phi\| ) $ tends in distribution to $N(0,\Sigma)$ as $n \rightarrow \infty$, where
\begin{align*}
\Sigma & = n \mbox{var} (\overline{Y}) = E(Y_1 Y_1^T) - E(Y_1) E(Y_1)^T \\
 &= d^{-1} \left[ \{ 1 - \eta_{2,d} (\|\phi\|) \} I + \{ (d+1) \eta_{2,d}(\| \phi \|) -1 -d \eta_{1,d}^2(\| \phi \|) \} \frac{\phi \phi^T}{\|\phi\|^2} \right]. 
\end{align*}
The monotonicity of $\eta_{1,d}$ implies that the Delta Method is applicable to $\overline{Y}$ and we have $\sqrt{n} (\hat{\phi}_{MM} - \phi) \stackrel{d}{\rightarrow} N(0,\Lambda^T \Sigma \Lambda) \ (n \rightarrow \infty)$, where
$$
\Lambda= \frac{\partial }{\partial \beta} \eta_{1,d}^{-1} (\|\beta \|) \frac{\beta}{\|\beta\|} \Biggr|_{\beta=\eta_{1,d} (\|\phi\|) \phi/\|\phi\|} = {\eta_{1,d}^{-1}}' \{ \eta_{1,d} (\|\phi\|) \} \frac{\phi^T \phi}{\|\phi\|^2} + \frac{\|\phi\|}{|\eta_{1,d} (\|\phi\|)|} \left( I - \frac{\phi^T \phi}{ \|\phi\|^2} \right). \label{eq:lambda}
$$
Here
\begin{equation}
{\eta_{1,d}^{-1}}' \{ \eta_{1,d} (\|\phi\|) \} = \frac{\partial}{\partial  x} \eta_{1,d}^{-1} (x) \Biggr|_{x=\eta_{1,d}(\|\phi\|)} = \frac{1}{\frac{\partial }{\partial y} \eta_{1,d}(y) } \Biggr|_{y=\|\phi\|}. \label{eq:eta1_dash_proof}
\end{equation}
The first derivative of $\eta_{1,d}$ is
$$
\frac{\partial }{\partial y} \eta_{1,d}(y) = \frac{\partial }{\partial y} E_Z \left( \frac{b-Z}{bZ-1} \right) = - b' E_Z \left\{ \frac{1-Z^2}{(bZ-1)^2} \right\},
$$
where $b=-2 y/(1+y^2)$, $b'=(\partial /\partial y) b = -2(1-y^2)/(1+y^2)^2$, and $Z$ follows the symmetric beta distribution (\ref{eq:beta}).
It follows from \citet[equation 9.111]{gra} that
\begin{align*}
E_Z \left\{ \frac{1-Z^2}{(bZ-1)^2} \right\} & = \int_{-1}^{1} \frac{1-z^2}{(bz-1)^2} \frac{(1-z^2)^{(d-2)/2}}{B(d/2,1/2)} dz \\
 &= \frac{d}{d+1} \frac{1}{(1+b)^2} F \left( 2,\frac{d}{2}+1 ; d+2 ; \frac{2b}{1+b} \right). 
\end{align*}
Substituting these results into (\ref{eq:eta1_dash_proof}), we obtain the expression for ${\eta_{1,d}^{-1}}' \{ \eta_{1,d} (\|\phi\|) \}$ given in Theorem \ref{thm:mme}.
\end{proof}

\subsection{Proof of Theorem \ref{thm:n3}}

\begin{proof}
\begin{enumerate}
\item[(i)] Clearly, the likelihood function is unbounded at $\phi=y_1$ and bounded otherwise.

\item[(ii)]
Let $y_1=e_1$ and $y_2=-e_1$ be the observations from $C^*_d(\phi)$.
Then the likelihood function $L_2$ is proportional to
\begin{align*}
L_2(\phi) & \propto \left[ \frac{(1-\|\phi\|^2)^2}{ \| e_1 - \phi \|^2 \| -e_1 - \phi\|^2  } \right]^d \\
 & = \left[ \frac{(1-\phi_1^2-\|\tilde{\phi}\|^2)^2}{ \{ (1-\phi_1)^2 + \|\tilde{\phi}\|^2 \} \{ (1+\phi_1)^2 + \|\tilde{\phi}\|^2 \} } \right]^d
\end{align*}
where $\phi=(\phi_1,\tilde{\phi}^T)^T,\ \tilde{\phi}=(\phi_2,\ldots,\phi_{d+1})^T$.
It follows from this expression that the maximum likelihood estimate of $\tilde{\phi}$ has to be zero.
Then the contour of the maximum likelihood of $\phi$ is given by the line
\begin{equation}
\sk{C} = \{ (\beta,0,\ldots,0)^T \,;\, \beta \in \mathbb{R} \}. \label{eq:l}
\end{equation}
Rotation of $e_1$ implies that, for $y_1=-y_2$, the contour of the maximum likelihood is the line connecting $y_1$ and $y_2$.

When $y_1=y_2$, it is clear from Theorem \ref{thm:n3}(i) that $\hat{\phi}_{ML}=y_1$.

Consider the the maximum likelihood for general $y_1$ and $y_2$ $(y_1 \neq \pm y_2)$.
Let $\check{\phi}=(0,\sk{\gamma},0,\ldots,0)^T$, where $\sk{\gamma}= ( [ 1-\{ (1-y_1^T y_2)/2 \}^{1/2} ]/[ 1 + \{(1-y_1^T y_2)/2 \}^{1/2}] )^{1/2} $.
Then ${\cal M}_{I,\check{\phi}} (e_1)^T {\cal M}_{I,\check{\phi}} (-e_1)= y_1^T y_2 $.
This implies that the angle between $y_1$ and $y_2$ is the same as that between ${\cal M}_{I,\check{\phi}} (e_1)$ and ${\cal M}_{I,\check{\phi}} (-e_1)$.
The contour of the maximum likelihood \sk{$C$} in (\ref{eq:l}) is transformed via ${\cal M}_{\check{\phi},I}$ onto
$$
\left\{(\phi_1,\phi_2,0,\ldots,0)^T \,;\, \phi_1^2 + \{\phi_2 - (1+\sk{\gamma}^2)/(2\sk{\gamma})\}^2 = \{(1-\sk{\gamma}^2)/(2\sk{\gamma})\}^2 \right\}.
$$
Therefore the contour of the maximum likelihood for general $y_1$ and $y_2$ $(y_1 \neq \pm y_2)$ is the circle perpendicular to the unit sphere with chord $(y_1,y_2)$ in the two-dimensional plane spanned by $y_1$ and $y_2$.

\item[(iii)]
In order to obtain the maximum likelihood estimator for the spherical Cauchy for $n=3$, we first consider the maximum likelihood estimator for the \sk{$d$-variate} $t$-distribution with $d$ degrees of freedom (\ref{eq:euclid_c2}) and then transform the maximum likelihood estimator using Theorem \ref{thm:mle}.

First we consider the observations $x_1=-e_1, \ x_2=0$ and $x_3=e_1$ sampled from the \sk{$d$-variate} $t$-family with $d$ degrees of freedom $C_d(\theta)$ given in (\ref{eq:euclid_c2}).
In a similar manner to \cite{fer} and \cite{mcc96}, it can be seen that the maximum likelihood estimate of $\theta$ is given by $\hat{\theta} = 0+ (1/\sqrt{3}) i$.

Next we transform $\hat{\theta}$ in order to obtain the maximum likelihood estimate of $\theta$ for general $x_1,x_2$ and $x_3$.
\citet{let} showed that the family (\ref{eq:euclid_c2}) is closed under the following transformation (\ref{eq:mobius_r}).
Let $\theta = \mu + i \sigma \in (\mathbb{R}^d + i \mathbb{R} ) \cup \{\infty\}$.
Define the following operations
$$
\theta+ a = \mu + a + i \sigma, \quad \gamma \theta = \gamma \mu + i \gamma \sigma, \quad 
$$
$$
\|\theta\| = \{ \|\mu\|^2 + \sigma^2 \}^{1/2}, \quad A \theta = A \mu + i \sigma,
$$
where $a, \gamma $ and $A$ are defined as in (\ref{eq:mobius_r}).
Then, if $X \sim C_d(\theta)$,
$$
X + a \sim C_d(\theta + a), \quad \gamma X \sim C_d (\gamma \theta), \quad A X \sim C_d (A \theta), \quad \frac{X}{\|X\|^2} \sim C_d \left( \frac{\theta}{\|\theta\|^2} \right). \label{eq:cauchy_properties}
$$

Here we transform $(-e_1,0,e_1)$ to the general $(x_1,x_2,x_3)$.
To achieve this, we first set $(\tilde{x}_1,\tilde{x}_3) = R (x_1-x_2,x_3-x_2)$, where $R$ is a $d \times d$ rotation matrix such that $ \tilde{x}_j = (\tilde{x}_{j1},\tilde{x}_{j2},0,\ldots,0)^T$ $(j=1,3)$.
Note that $(x_1,x_2,x_3)$ and $(\tilde{x}_1,0,\tilde{x}_3)$ constitute the same triangle apart from translation and rotation.
Then the three points $(-e_1,0,e_1)$ are transformed to $(\tilde{x}_1,0,\tilde{x}_3)$ via the transformation
\begin{equation}
{\cal E}(t) = A \left( \gamma \, \frac{t+a}{\|t+a\|^2} + b \right), \quad t \in \mathbb{R}^d \times i \mathbb{R}, \label{eq:f}
\end{equation}
where
$$
A = \left(
\begin{array}{cc}
\tilde{A} & O \\
O & I
\end{array}
\right),
\quad
\tilde{A} =  \frac{1}{|\alpha \beta|} \left(
\begin{array}{cc}
- \mbox{Re}(\alpha \beta) & - \mbox{Im} (\alpha \beta) \\
- \mbox{Im} (\alpha \beta) & \mbox{Re}(\alpha \beta)
\end{array}
\right),
$$
$$
a = \left( \mbox{Re}( \beta), \mbox{Im} ( \beta), 0, \ldots, 0 \right)^T, \quad b= - \left| \alpha / \beta \right| \left( \mbox{Re}( \beta ), \mbox{Im} ( \beta), 0, \ldots, 0 \right)^T,
$$
$$
\gamma = |\alpha \beta|, \quad \alpha = \frac{2 z_1 z_3}{z_1 + z_3}, \quad \beta = \frac{z_1-z_3}{z_1+z_3}, \quad z_j = \tilde{x}_{j1} + i \tilde{x}_{j2}, \quad j=1,3.
$$
Then $x_1$ and $x_3$ have the expression
$$
x_1 = x_2 + R^T {\cal E}(-e_1) \quad \mbox{and} \quad x_3 = x_2 + R^T {\cal E}(e_1).
$$
Substituting $t=0 + (1/\sqrt{3}) i$ into (\ref{eq:f}), the estimate of the parameter for general $(x_1,x_2,x_3)$ is given by $\hat{\theta} = x_2 + R^T {\cal E}(t)$.
After some algebra, it can be seen that the estimate is of the form
$$
\hat{\mu} = \frac{\|x_1-x_2\|^2 x_3 + \|x_2-x_3\|^2 x_1 + \|x_3-x_2\|^2 x_1}{ \|x_1-x_2\|^2 + \|x_2-x_3\|^2 + \|x_3-x_1\|^2 },
$$
$$
\hat{\sigma} = \pm \sqrt{3} \frac{\|x_1-x_2\| \|x_2 - x_3\| \|x_3-x_1\|}{\| x_1 - x_2\|^2 + \|x_2 - x_3\|^2 + \|x_3 - x_1\|^2 }.
$$
Finally the maximum likelihood estimate of $\phi$ for the sample $y_1,y_2$ and $y_3$ can be obtained via Theorem \ref{thm:mle} by substituting $x_j$ into ${\cal P} (y_j)$ and transforming $\hat{\mu} + i \hat{\sigma}$ via the mapping ${\cal P}^{-1}$.
\end{enumerate}
\end{proof}

\subsection{Proof of Lemma \ref{pro:fi}}

\begin{proof}
The Fisher information matrix is expressed as
\begin{equation}
\begin{split} \label{eq:fi}
{\cal I} = & \ 2d \left\{ \frac{1}{1-\|\phi\|^2} + E \left( \frac{1}{\|Y-\phi\|^2} \right) \right\} I \\
 & + 4d \left[ \frac{\phi \phi^T}{(1-\|\phi\|^2)^2} - E \left\{ \frac{(Y-\phi) (Y-\phi)^T}{\|Y-\phi\|^4} \right\} \right].
\end{split}
\end{equation}
Put $Z={\cal M}_{I,-\phi}(Y)$.
Then it follows from Theorem \ref{thm:closure_cauchy_s} that $Z$ has the uniform distribution on the sphere.
Therefore
$$
E_Y \left( \frac{1}{\|Y-\phi\|^2} \right) = \frac{E_Z (\|Z+\phi \|^2)}{(1-\|\phi\|^2)^2} = \frac{1+\|\phi\|^2}{(1-\|\phi\|^2)^2}
$$
and
$$
E_Y \left\{ \frac{(Y-\phi) (Y-\phi)^T}{\|Y-\phi\|^4} \right\} = \frac{E_Z \left\{ (Z+\phi) (Z+\phi)^T \right\}}{(1-\|\phi\|^2)^2}  = \frac{1}{(1-\|\phi\|^2)^2} \left\{ \frac{1}{d+1} I + \phi \phi^T \right\}.
$$
Substituting these results into (\ref{eq:fi}), we obtain Lemma \ref{pro:fi}.
\end{proof}

\subsection{Proof of Theorem \ref{thm:unimodality}}

\begin{proof} 
Let $\phi^*$ be a stationary point of the loglikelihood function (\ref{eq:loglikelihood}).
For convenience, write
\begin{equation}
z_j = {\cal M}_{I,-\phi^*} (y_j). \label{eq:z_j}
\end{equation}
It holds that $z_j \in S^d$.
Then the estimating equation for $\phi$ can be simply expressed as
\begin{equation}
\sum_{j=1}^n z_j = 0. \label{eq:z_sum}
\end{equation}
The second derivative of the loglikelihood function is
\begin{align*}
\frac{\partial^2 \ell }{\partial \phi \partial \phi^T} = -2d \left\{ \frac{n I }{1-\|\phi\|^2} + \frac{2n \phi \phi^T }{(1-\|\phi\|^2)^2}  + \sum_{j=1}^n \frac{I}{\|y_j-\phi \|^2}  - 2 \sum_{j=1}^{n} \frac{(y_j-\phi) (y_j-\phi)^T}{ \| y_j-\phi \|^4}   \right\}.
\end{align*}
Using $z_j$ defined in (\ref{eq:z_j}), the second derivative of the loglikelihood function at $\phi = \phi^*$ can be expressed as
\begin{align*}
\lefteqn{ \frac{\partial^2 \ell}{\partial \phi \partial \phi^T} \biggr|_{\phi=\phi^*} } \nonumber \\
 & = -2d \left[ \left\{ \frac{n}{1-\|\phi^*\|^2} + \sum_{j=1}^n \frac{\|z_j+\phi^*\|^2}{(1-\|\phi^*\|^2)^2} \right\} I + \frac{2n \phi^* \phi^{*T}}{(1-\|\phi^* \|^2)^2} -2 \sum_{j=1}^n \frac{(z_j+\phi^*)(z_j+\phi^*)^T}{(1-\|\phi^* \|^2)^2} \right] \nonumber \\
 & = - \frac{4d}{(1-\|\phi^* \|^2)} \left( n I - \sum_{j=1}^n z_j z_j^T \right). 
\end{align*}
The second equality follows from equation (\ref{eq:z_sum}).
Therefore, for any $t \in S^d$,
$$
t^T \frac{\partial^2 \ell }{\partial \phi \partial \phi^T} \biggr|_{\phi=\phi^*} t  = - \frac{4d}{(1-\|\phi^* \|^2)} \left\{ n - \sum_{j=1}^n (z_j^T t)^2 \right\} <  - \frac{4d}{(1-\|\phi^* \|^2)} \left( n - n \right) = 0,
$$
implying the second derivative of the loglikelihood function at $\phi=\phi^*$ is a negative definite matrix.
Here the inequality follows from the assumption that $y_j \neq y_k$ for some $(j,k)$.
Hence any stationary point of the loglikelihood function is a local maximum.
\end{proof}

\subsection{Proof of Lemma \ref{lem:marginal}}

\begin{proof}
It follows from the last sentence in \S \ref{sec:marginal} that $\eta_{k,\nu}(\rho)$, the $k$th moment of the \sk{random variable having the} distribution (\ref{eq:marginal}), can be expressed as
\begin{equation}
\eta_{k,\nu}(\rho) = E_{Z} \left\{ \left( \frac{b-Z}{b Z -1}  \right)^k \right\}, \label{eq:mu_k_alternative}
\end{equation}
where $b =-2 \rho/(1+\rho^2)$ and $Z$ follows the symmetric beta distribution (\ref{eq:beta}).
With this expression,
$$
\frac{\partial}{\partial \rho} \eta_{k,\nu}(\rho) = E_Z \left\{ \frac{\partial}{\partial \rho} \left( \frac{b -Z}{b Z -1}  \right)^k  \right\} = 2 k \frac{1-\rho^2}{(1+\rho^2)^2} E_Z \left\{ \frac{(b-Z)^{k-1} (1-Z^2)}{(bZ-1)^{k+1}} \right\} >0.
$$

The symmetry of the distribution of $Z$ implies that $\eta_{k,\nu}(0)=0$.
As for the limit of $\eta_{k,\nu}(\rho)$, we see that, for fixed $z \, (\in (-1,1))$, $\lim_{\rho \rightarrow 1} \{(b-z)/(b z -1)\} =1$.
Therefore this fact and the dominated convergence theorem imply that $\lim_{\rho \rightarrow 1} \eta_{k,\nu}(\rho) = \lim_{\rho \rightarrow 1} E_Z [ \{(b-z)/(b z -1)\}^k ] = 1$.
\end{proof}

\subsection{Proof of Theorem \ref{thm:1st_moment}}

\begin{proof}
\begin{enumerate}
\item[(i)] It follows from (\ref{eq:mu_k_alternative}) that the mean of $Y_1$ can be expressed as
\begin{equation}
\eta_{1,\nu} (\rho) = E_{Y_1} (Y_1) = E_{Z} \left( \frac{b-Z}{b Z -1}  \right) = -  \frac{1}{b} -  \frac{1-b^2}{b^2} E_{Z} \left( \frac{1}{Z-b^{-1}} \right), \label{eq:mu_1_alternative}
\end{equation}
where $b=-2 \rho/(1+\rho^2)$ and $Z$ follows the symmetric beta distribution (\ref{eq:beta}).
With this representation, the first equality in Theorem \ref{thm:1st_moment}(i) \sk{is} clear from equation (9.111) of \cite{gra}.
The second equality in Theorem \ref{thm:1st_moment}(i) follows from equations (9.131.1) and (9.134.1) of \cite{gra}.

\item[(ii)] 
The closed-form expressions for $\eta_{1,1} (\rho)$ and $\eta_{1,2} (\rho)$ are available by applying equations (9.131.1) and (9.121.6) of \cite{gra}, respectively, to Theorem \ref{thm:1st_moment}(i).
The integral representation of the hypergeometric series \citep[equation 9.111]{gra} leads to closed-form expressions for $\eta_{1,3} (\rho)$ and $\eta_{1,4} (\rho)$.

\item[(iii)] 
Gauss recursion formulas (9.137.5) and (9.137.15) of \cite{gra} imply that
\begin{align*}
F \left( \frac12, \frac{\nu-1}{2}; \frac{\nu+1}{2}; z \right) = \ & \frac{\nu-1}{(\nu-2) (\nu-3)} \biggl\{ \left( \nu -2 - \frac{\nu-3}{z} \right) F \left( \frac12, \frac{\nu-3}{2}; \frac{\nu-1}{2}; z \right)  \\
& + \frac{\nu-3}{z} F \left(\frac12, \frac{\nu-5}{2}; \frac{\nu-3}{2}; z \right)  \biggr\}.
\end{align*}
The recursion formula presented in Theorem \sk{\ref{thm:1st_moment}(iii)} can be obtained \sk{by substituting the equation above into the second expression of $\eta_{1,\nu}(\rho)$ in} Theorem \ref{thm:1st_moment}(i).
\end{enumerate}
\end{proof}

\subsection{Proof of Theorem \ref{thm:2nd_moment}}

\begin{proof}
\begin{enumerate}
\item[(i)]
It follows from equation (9.111) of \citet{gra} and equation (\ref{eq:mu_1_alternative}) that $\eta_{2,\nu} (\rho)$ has the form (\ref{eq:2nd_moment}).

\item[(ii)] 
For $\nu=1,\ldots,4$, Theorem \ref{thm:1st_moment} \sk{implies} that the hypergeometric series $F \{ 1 , \nu/2 ; \nu; - 4 \rho/(1-\rho)^2 \}$ in the second term in the right-hand side of (\ref{eq:2nd_moment}) can be expressed in \sk{closed} form.
The hypergeometric series $F \{ 2 , \nu/2 ; \nu ; - 4 \rho/ (1-\rho)^2 \}$ in the third term of the right-hand side of (\ref{eq:2nd_moment}) can be calculated partly using its integral representation \citep[equation 9.111]{gra}.
Summarizing these facts, the closed-form expressions for $\eta_{2,\nu} (\rho)$ are available for $\nu=1,\ldots,4$.

\item[(iii)]
From equation (9.134.3) of \cite{gra}, we have
$$
F \left\{ 2 , \frac{\nu}{2}, \nu, -\frac{4\rho}{(1-\rho)^2} \right\} = \frac{(1-\rho)^2}{(1+\rho)^2} F \left( 1 , \frac{\nu-2}{2} , \frac{\nu+1}{2} , z \right),
$$
\sk{where $z=-4 \rho^2 / (1-\rho^2)^2$.}
Then it follows from this equation and Theorem \ref{thm:2nd_moment}(i) that
\begin{equation}
F \left( 1 , \frac{\nu-2}{2} , \frac{\nu+1}{2} , z \right) = \frac{(1+\rho^2)^2}{(1-\rho^2)^2} \left\{ 1 - \frac{4 \rho}{1+\rho^2} \eta_{1,\nu} (\rho) + \frac{4 \rho^2}{(1+\rho^2)^2} \eta_{2,\nu} (\rho) \right\}. \label{eq:f1}
\end{equation}
Equations (9.137.5) and (9.137.15) of \cite{gra} imply that
\begin{align}
\begin{split} \label{eq:f1_recursion}
F \left( 1 , \frac{\nu-2}{2} , \frac{\nu+1}{2} , z \right) = \ & \frac{1}{(\nu-3)(\nu-4)} \biggl\{ c_1 F \left( 1 , \frac{\nu-4}{2} , \frac{\nu-1}{2} , z \right) \\
& - c_2 F \left( 1 , \frac{\nu-6}{2} , \frac{\nu-3}{2} , z \right) \biggr\}.
\end{split}
\end{align}
Substituting (\ref{eq:f1}) into (\ref{eq:f1_recursion}), we obtain Theorem \ref{thm:2nd_moment}(iii).
\end{enumerate}
\end{proof}

\end{document}